\documentclass[mathpazo]{cicp}

\usepackage[table,xcdraw]{xcolor}
\usepackage{float}
\usepackage{caption}
\usepackage{subfigure}
\usepackage{graphicx}
\usepackage{float}
\usepackage{comment}
\usepackage{multirow}
\usepackage{mathrsfs}
\usepackage{amsmath,amssymb} 
\usepackage{color}
\usepackage{times}
\usepackage{epsfig}
\usepackage{diagbox}
\usepackage{afterpage}
\usepackage{makecell}
\usepackage{lineno}
\usepackage{amsthm}   

\usepackage{bm}
\usepackage{booktabs}       
\usepackage{hyperref}
\hypersetup{colorlinks=true,linkcolor=black}
\usepackage{docmute}

\usepackage{subfiles}  
\usepackage{tikz}
\usepackage{calc}
\usepackage[ruled,linesnumbered]{algorithm2e}
\usepackage{CJK}

\newcommand{\be}{\begin{eqnarray}}
\newcommand{\ee}{\end{eqnarray}}
\newcommand{\beq}{\begin{equation}\begin{aligned}}
\newcommand{\eeq}{\end{aligned}\end{equation}}
\newcommand{\beqn}{\begin{equation*}\begin{aligned}}
\newcommand{\eeqn}{\end{aligned}\end{equation*}}
\newcommand{\ben}{\begin{eqnarray*}}
\newcommand{\een}{\end{eqnarray*}}

\newtheorem{thm}{Theorem}
\newcommand{\paramA}{1}
\newcommand{\paramNumBlock}{10}

\graphicspath{{D:/research/PDE+DL/like/codes/AdaptiveDNN/AdaptiveDNN/results_block\paramNumBlock/results_cluster_a_\paramA/figures2/}}
\begin{document}
\title{Adaptive Basis-inspired Deep Neural Network for Solving Partial Differential Equations with Localized Features}
\author[Li K et.~al.]{Ke Li\affil{1},
	Yaqin Zhang\affil{2}\comma\affil{3},
	Yunqing Huang\affil{2}, 
	Chenyue Xie\affil{4}\comma\footnotemark[1], and Xueshuang Xiang\affil{3}\comma\corrauth}
\address{
	\affilnum{1}\ Information Engineering University, Zhengzhou 450001, P. R. China.\\
	\affilnum{2}\ School of Mathematics and Computational Science,
	Xiangtan University,
	Xiangtan 411105, P. R. China. \\
	\affilnum{3}\ Qian Xuesen Laboratory of Space Technology,
	China Academy of Space Technology, Beijing 100094, P. R. China.\\
	\affilnum{4}\ Department of Modern Mechanics, University of Science and Technology of China, Hefei, 230002, P. R. China.\\}

\emails{
	{\tt cyxie@ustc.edu.cn} (C.~Xie),
	{\tt xiangxueshuang2023@163.com} (X.~Xiang)}
\begin{abstract}
This paper proposes an Adaptive Basis-inspired Deep Neural Network (ABI-DNN) for solving partial differential equations with localized phenomena such as sharp gradients and singularities. Like the adaptive finite element method, ABI-DNN incorporates an iteration of “solve, estimate, mark, enhancement”, which automatically identifies challenging regions and adds new neurons to enhance its capability. A key challenge is to force new neurons to focus on identified regions with limited understanding of their roles in approximation. To address this, we draw inspiration from the finite element basis function and construct the novel Basis-inspired Block (BI-block), to help understand the contribution of each block. With the help of the BI-block and the famous Kolmogorov Superposition Theorem, we first develop a novel fixed network architecture named the Basis-inspired Deep Neural Network (BI-DNN), and then integrate it into the aforementioned adaptive framework to propose the ABI-DNN. Extensive numerical experiments demonstrate that both BI-DNN and ABI-DNN can effectively capture the challenging singularities in target functions. Compared to PINN, BI-DNN attains significantly lower relative errors with a similar number of trainable parameters. When a specified tolerance is set, ABI-DNN can adaptively learn an appropriate architecture that achieves an error comparable to that of BI-DNN with the same structure.
\end{abstract}
\maketitle

\section{Introduction}
Recently, the remarkable success of Deep Neural Networks (DNNs) in data science has motivated numerous researchers to explore their potential in solving Partial Differential Equations (PDEs). Various DNN-based methods have emerged, such as the deep Ritz method (DRM) \cite{weinan2018deep}, the Physical Information Neural Network (PINN) \cite{PINN}, and the deep Galerkin method (DGM) \cite{Sirignano_Spiliopoulos_2018}.
In particular, the PINN introduced by Raissi et al. has attracted much attention due to its ability to embed physical equations into the network architecture, allowing for accurate predictions without the requirement of extensive data. It offers advantages over traditional mesh-based methods, including flexibility with complex geometries, efficiency in handling high-dimensional problems, and ease of implementation. As a result, the PINN has become an increasingly important tool for scientists and engineers seeking to understand better and predict the behavior of complex systems, such as fluid mechanics\cite{Cai_Mao_Wang_Yin_Karniadakis_2021, Xie_Wang_Li_Wan_Chen_2019, Xie_Wang_E_2020}, hydrogeophysics\cite{Meng_Karniadakis_2019} and transport phenomena in porous media\cite{Pang_Lu_Karniadakis_2018}.

Despite their successes, PINNs still face challenges in dealing with complicated PDE problems. For instance, the training of PINNs often encounters difficulties when target solutions involve sharp gradients or discontinuities. 
Therefore, continuous refinements are conducted from diverse perspectives to further improve the performance of PINNs. For instance, \cite{Wu_Zhu_Tan_Kartha_Lu_2022, Gao_Yan_Zhou_2022, Gao_Tang_Yan_Zhou_2023} proposed various adaptive sampling strategies that automatically adjust collocation points according to the residual of the PDE, the gradient information of the NN or insights derived from them. \cite{Jagtap_Kawaguchi_Karniadakis_2020, Jagtap_Kawaguchi_Em_Karniadakis_2020} introduced adaptive activation functions with scalable parameters to enhance convergence rates and solution accuracy.  
\cite{McClenny_Braga-Neto_2020, Xiang_Peng_Liu_Yao_2022} focused on adaptive loss weighting strategies 
to achieve balanced training. 
\cite{Taylor_Bastidas_Calo_Pardo_2024} addressed problems with singular solutions by adaptive domain decomposition.

Besides factors explored in the above studies, the choice of neural network architecture, including the number of layers, the number of neurons within each layer, and patterns of connectivity, is also crucial for an accurate scheme \cite{Cai_Chen_Liu_2022}. However, network architecture design is often empirical due to a limited understanding of the DNN. 
Noticing that overly simple architectures might underperform and complex ones might overfit with limited data, designing an appropriately complex architecture is inherently challenging.
An attractive idea is to shift from a pre-specified architecture with fixed complexity to an adaptively generated, task-specific network architecture, which is more efficient and suitable for the challenges it addresses.
Regrettably, literature on the adaptive network architecture design within PINNs is still limited. A pioneering effort is the Adaptive Network Enhancement (ANE) method\cite{Liu_Cai_2022, Liu_Cai_Chen_2022, Cai_Chen_Liu_2022}, which starts with a small network and iteratively adds neurons based on the a posteriori estimator until the desired accuracy is reached. 
However, the ANE method heavily relies on the physical partition generated by the ReLU activation function, and computing the physical partition for a deep network is computationally intensive.

In this paper, we focus on the adaptive design of network architecture and draw inspiration from the adaptive finite element method (AFEM) to propose the Adaptive Basis-inspired Deep Neural Network (ABI-DNN) for solving problems with localized features. The advantage of AFEMs lies in their ability to locally refine the mesh, which enhances the approximating capacity of the finite element function space, so that better resolution of local features can be obtained. 
Similarly, with the function space underlying a neural network determined by its architecture, incrementally augmenting the network with new neurons is promising to obtain a more precise adaptation to the characteristics of the target function. 
Unlike AFEM, which employs the a posteriori estimator to directly guide mesh refinement, it is not straightforward to determine how to direct newly added neurons to focus on the marked high-error regions. This challenge arises from the limited understanding of the role each neuron plays in the approximation process.
To address this difficulty, we start from the one-dimensional case and introduce the novel Basis-inspired Block (BI-block), which is designed to emulate the properties of linear finite element basis functions. Since the attention region of each BI-block with proper initialization is well understood, BI-blocks that focus on local regions with high errors can be dynamically integrated to enhance the network's approximating ability for a more precise capture of the localized characteristics.
To avoid tensor products in high-dimension problems, we draw inspiration from the famous Kolmogorov Superposition Theorem (KST),  which allows the representation of any continuous multivariate function as sums and compositions of simpler univariate functions, offering a powerful framework for constructing complex multivariate representations using univariate components.  Following the framework of KST, we construct the proposed Basis-inspired DNN (BI-DNN), with BI-blocks as the foundational building blocks.
With this innovative architecture established, we develop an adaptive enhancement strategy and propose the adaptive basis-inspired deep neural network (ABI-DNN). Similar to the adaptive finite element method, our ABI-DNN initiates with an initial small architecture that yields coarse solutions upon training. It then iteratively and dynamically integrates new BI-blocks into the architecture guided by the error distribution given by the current approximating solution.  After that, we retrain the augmented model until the prescribed tolerance is attained.
Various numerical experiments are conducted to validate the performance of the proposed ABI-DNN. In summary, the main contributions of this work are as follows:
\begin{itemize}
	\item We propose a novel network architecture, BI-DNN, which is inspired by the FEM basis functions and the KST, demonstrating superior performance compared to PINNs in addressing problems characterized by localized features.
	\item We introduce an adaptive network architecture enhancement framework, ABI-DNN, including practical strategies for solving, estimation, marking, and enhancement. This framework enables the network to adaptively scale from a small-sized network to meet the predefined tolerance and achieve improved resolution in challenging regions.
	\item We successfully address a series of challenging problems, including those with singularities or high frequencies. The numerical results demonstrate that the BI-DNN achieves errors approximately one or two orders of magnitude lower than the standard PINN (see Figs. \ref{fig:fitting1d:params}, \ref{fig:pde2d:params:onepeak}, and  \ref{fig:pde2d:params:twopeaks}). Furthermore, the ABI-DNN can adaptively generate a suitable network architecture to meet the prescribed tolerance and typically achieves an error comparable to or smaller than that of the BI-DNN with the same architecture (see Tables \ref{tab:fitting1d:ABI-DNN:singular}, \ref{tab:fitting1d:ABI-DNN:smooth}, \ref{tab:pde2d:ABI-DNN:onepeak}, \ref{tab:pde2d:ABI-DNN:Lshape}, and  \ref{tab:pde2d:ABI-DNN:burgers}).
\end{itemize}

The remainder of this paper is organized as follows. Section 2 provides a concise overview of the standard PINN. Section 3 introduces the novel BI-block and BI-DNN. Section 4 presents the ABI-DNN method, which is an adaptive network architecture enhancement framework. Section 5 details extensive numerical experiments to demonstrate the efficiency of the proposed method. Finally, Section 6 concludes the paper with a brief summary.

\section{Physics-informed Neural Network}
Physics-informed neural networks\cite{PINN} are DNN-based methods for solving PDEs as well as inverse problems of PDEs, which have received much attention recently. Specifically, let's consider a general Partial Differential Equation (PDE)  represented as:
\beq \label{eq:general_equation}	
\mathcal{L} (u) &= f, \quad \text{in}\  \Omega, \\
\mathcal{B} (u) &= g, \quad \text{on}\  \partial\Omega, 
\eeq
where $\mathcal{L}$ is the differential operator such as $-\Delta$, and $\mathcal{B}$ denotes the boundary operator that imposes conditions like Dirichlet or Neumann conditions. We remark that when addressing time-dependent problems, time $t$ can be regarded as an additional coordinate in $x$, and $\Omega$ denotes the spatio-temporal domain. Consequently, the initial condition can be simply regarded as a special type of Dirichlet boundary condition.
In standard PINN, the target function $u(x)$ is approximated by a feedforward neural network $u(x;\theta)$ with parameters $\theta$. In general, a feedforward neural network with $L$ hidden layers can be mathematically represented by a composite function of the form: 
\begin{equation}
	u(x;\theta) = {\bf F}_{L+1} \circ \sigma \circ {\bf F}_{L} \circ \sigma \circ\cdots \circ {\bf F}_{2} \circ \sigma \circ {\bf F}_{1}(x), \label{eq:DNN}
\end{equation}
where $\sigma$ denotes an activation function such as the hyperbolic tangent (Tanh) or rectified linear unit (ReLU), and each ${\bf F}_{l},1\leq l \leq L$ represents an affine transformation defined as: 
\begin{equation*}
	{\bf F}_{l}(x)=W_l x+b_l.
\end{equation*}
Here, $W_l \in \mathbb{R}^{n_{l} \times n_{l-1}}, b \in \mathbb{R}^{n_{l}}$, $n_l$ is the number of neurons in the $l$-th layer, and $\theta:=\{W_l,b_l\}$ is the collection of all trainable parameters. The objective of PINN in addressing problem \eqref{eq:general_equation} is to find the optimal parameters 
$\theta$ that minimize the following discrete loss function: 
\begin{equation}
	\theta^{*} = argmin_{\theta} 
	\mathcal{M}_{\Omega}(\theta) + \beta * \mathcal{M}_{\partial \Omega}(\theta). \label{eq:loss}
\end{equation}
In this formulation, the discrete residual terms $\mathcal{M}_{\Omega}(\theta)$ and $\mathcal{M}_{\partial \Omega}(\theta)$  are defined as: 
\[
\mathcal{M}_{\Omega}(\theta):=\frac{1}{N_{r}} \sum_{i=1}^{N_{r}}|\mathcal{L} (u(x_{r}^{i};\theta))-f(x_{r}^{i})|^{2},\quad
\mathcal{M}_{\partial \Omega}(\theta):=\frac{1}{N_{b}} \sum_{i=1}^{N_{b}}|\mathcal{B} (u(x_{b}^{i};\theta))-g(x_{b}^{i})|^{2},
\]
which are responsible for enforcing the PDE constraints and the boundary conditions, respectively. The hyperparameter $\beta$ is a weight that balances these two terms. The sets $\big \{x_{r}^{i}\big \}_{i=1}^{N_{r}}$ and $\big \{x_{b}^{i}\big \}_{i=1}^{N_{b}}$ denote the collocation points located inside the domain $\Omega$ and on its boundary $\partial \Omega$, respectively. Stochastic optimization algorithms such as Adam\cite{kingma2014adam} can be employed to optimize problem \eqref{eq:loss}.

\section{Basis-inspired Deep Neural Network}
In this section, we will introduce the novel architecture BI-DNN, which will be employed in the adaptive framework in the next section. This architecture incorporates specially designed BI-blocks, which are crafted to either precisely or approximately represent finite element basis functions, enabling the network to capture complex patterns and enhance its representational power.
\subsection{Basis-inspired Block}

\subsubsection{Basis-inspired Block with ReLU Activation Function}
\begin{figure}[htbp]
	\begin{center}
		\begin{tabular}{cc}
			\centering
			\quad\includegraphics[width=0.45\textwidth]{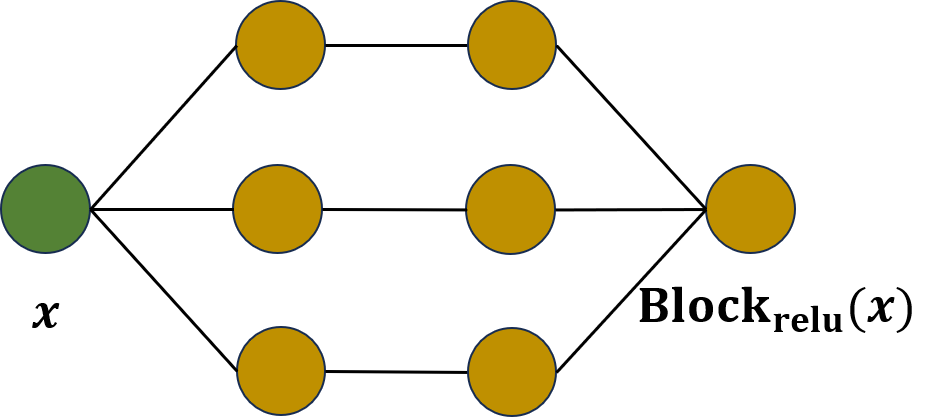}&
			\quad\includegraphics[width=0.45\textwidth]{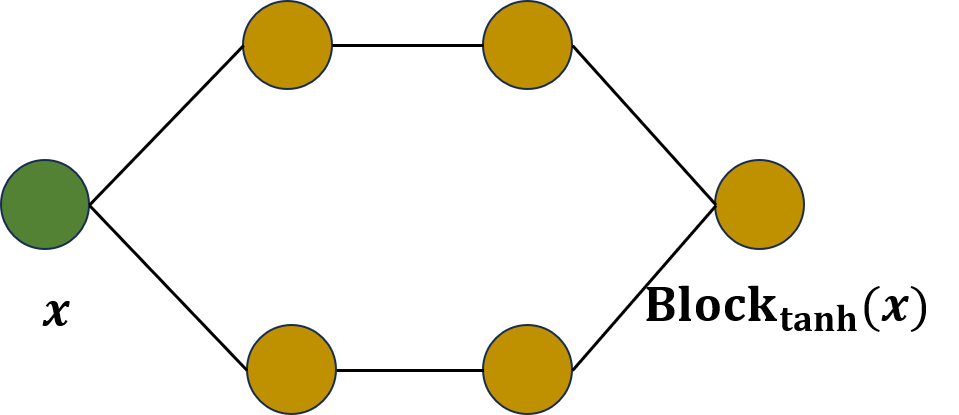}\\
			\hspace{-0.1in}
			(a) ReLU Basis Block\qquad & (b) Tanh Basis Block\qquad \\
		\end{tabular}
		\caption{
			Illustration of basis-inspired blocks with activation functions: (a) ReLU (Equation \eqref{eq:bi-block:relu:v2}) and (b) Tanh (Equation \eqref{eq:bi-blocks:tanh:v2}), respectively. }
		\label{fig:bi-blocks}
	\end{center}
\end{figure}

We start with the one-dimensional (1D) piecewise linear basis function to illustrate the construction of BI-blocks. 
In the conventional FEM, the computational domain is discretized into elements defined by a set of nodal points denoted as $x_j$ for $j=1,\cdots,n$. Associated with each nodal point is a basis function, enabling the representation of every finite element function as a linear combination of these nodal basis functions. The basis function associated with node $x_j$ can be written as: 
\ben
\varphi_j(x)=\left\{
\begin{array}{ll}
	\frac{x-x_{j-1}}{h_{j-1}}, & x_{j-1}\leq x< x_j,\\
	\frac{x_{j+1}-x}{h_j}, & x_j\leq x< x_{j+1},\\
	0, & otherwise,
\end{array}
\right.
\een
where $h_j=x_{j+1}-x_j$ is length of the $j$-th subinterval. 
In the domain of DNNs, these basis functions can be reinterpreted from a neural network perspective. In detail, a basis function can be elegantly reformulated using the ReLU activation, as detailed below:
\ben
\varphi_j(x)= \text{ReLU}(\frac{1}{h_{j-1}}(x-x_{j-1}))-2\text{ReLU}(\frac{1}{2}(\frac{1}{h_{j-1}}+\frac{1}{h_{j}})(x-x_j))+\text{ReLU}(\frac{1}{h_j}(x-x_{j+1})).
\een
This representation can further be compactly expressed as: 
\be
\varphi_j(x)=W_2\cdot\text{ReLU}(W_1\cdot x + b_1),
\label{eq:bi-block:relu:v1}
\ee
with the weights and biases defined as:
\ben
W_1=\left[\frac{1}{h_{j-1}},\frac{1}{2}(\frac{1}{h_{j-1}}+\frac{1}{h_{j}}),\frac{1}{h_{j}}\right]^T,\;
b_1=-\left[\frac{1}{h_{j-1}}x_{j-1},\frac{1}{2}(\frac{1}{h_{j-1}}+\frac{1}{h_{j}})x_j,\frac{1}{h_{j}}x_{j+1}\right]^T,\;
W_2=\left[1,-2,1\right].
\een
In essence, this formulation represents a neural network block consisting of a single hidden layer followed by a linear output layer, in which the weights and biases are the functions of nodal positions. Note that $W_1$ in \eqref{eq:bi-block:relu:v1} might be considerably large when a narrow basis function is necessary to capture sharp corners in the target functions. 
To alleviate potential training challenges associated with this, we rewrite the equation \eqref{eq:bi-block:relu:v1} as:
\be
\varphi_j(x) 
=W_2\cdot \text{ReLU}(W_1^2\cdot (W_1^1\cdot x + b_1^1)),
\label{eq:bi-block:relu:v2}
\ee
where the weights and the bias are defined as:
\be
W_1^1=\left[\sqrt{\frac{1}{h_{j-1}}},\sqrt{\frac{1}{2}(\frac{1}{h_{j-1}}+\frac{1}{h_{j}})},\sqrt{\frac{1}{h_{j}}}\;\right]^T,\quad
W_1^2=\begin{bmatrix}
	\sqrt{\frac{1}{h_{j-1}}} & 0 & 0\\
	0 &\sqrt{\frac{1}{2}(\frac{1}{h_{j-1}}+\frac{1}{h_{j}})} & 0\\
	0 & 0 & \sqrt{\frac{1}{h_{j}}}
\end{bmatrix},
\label{eq:bi-block:relu:weight:v2}
\ee
and 
\be
b_1^1=-\left[\sqrt{\frac{1}{h_{j-1}}}x_{j-1},\sqrt{\frac{1}{2}(\frac{1}{h_{j-1}}+\frac{1}{h_{j}})}x_j,\sqrt{\frac{1}{h_{j}}}x_{j+1}\right]^T.
\label{eq:bi-block:relu:bias:v2}
\ee
In the above formulation, We decompose the parameter $W_1$
into the square of its square root,  thereby effectively splitting one hidden layer into two successive layers.
This reformulation aims to alleviate
the training challenge associated with $W_1$ and allows for a more adaptive response to the geometric characteristics of the target functions.
For convenience, we refer to the neural network representation of the basis function depicted in \eqref{eq:bi-block:relu:v2}  as a `Basis-inspired Block' (abbreviated as BI-block)  with ReLU activation, denoted as $\text{Block}_{\text{relu}}(x;x_j,h_{j-1},h_j)$. We remark that the weights and biases within block $\text{Block}_{\text{relu}}(x;x_j,h_{j-1},h_j)$ are initialized using $x_j$, $h_{j-1}$, and $h_j$ according to equations  \eqref{eq:bi-block:relu:weight:v2} and \eqref{eq:bi-block:relu:bias:v2}, and they are allowed to evolve during the training process unless otherwise stated.
For illustration, a BI-block with ReLU activation is shown in Fig. \ref{fig:bi-blocks}(a).

\subsubsection{Basis-inspired Block with Tanh Activation Function}
Despite its popularity for simplicity, the ReLU activation function is not smooth enough for PINNs to be well-behaved. Critically, the backpropagation in PINNs involves calculating second-order derivatives of the neural network. However, the ReLU activation function, being piecewise linear, yields zero second derivatives except at a single point, which impedes the learning process\cite{Maczuga_Paszynski_2023}. Given the Tanh activation function's sufficient smoothness and effectiveness within the PINN framework, we will now turn our attention to employing it in the construction of the BI-block.

\begin{figure}[htbp]
	\begin{center}
		\begin{tabular}{cc}
			\centering
			\includegraphics[width=0.45\textwidth]{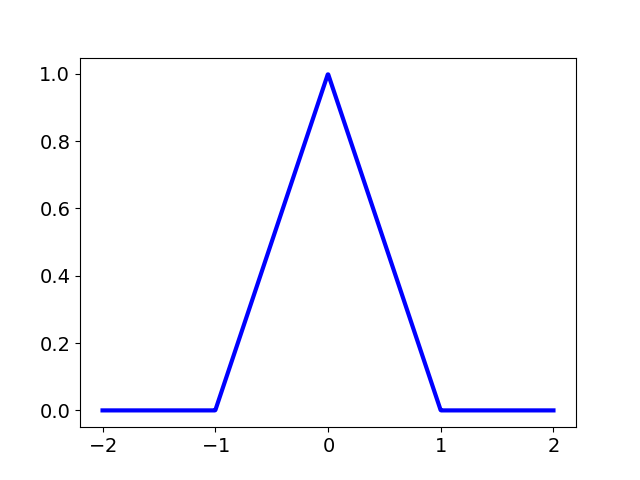}&
			\includegraphics[width=0.45\textwidth]{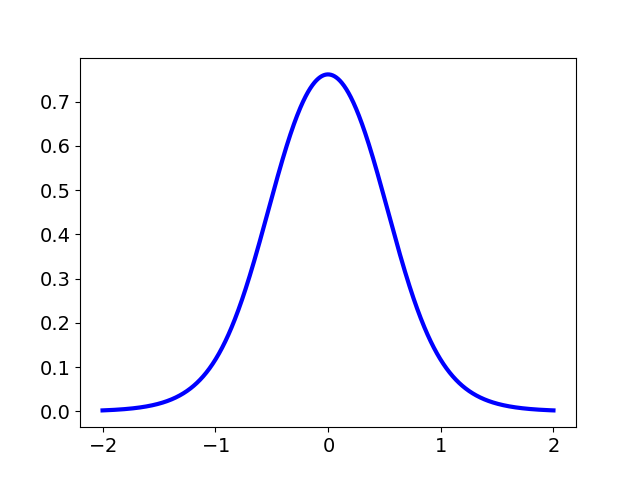}\\
			(a) & (b) \\
		\end{tabular}
		\caption{Image of basis-inspired blocks with (a) ReLU and (b) Tanh activation functions, respectively. Here, equally spaced nodes within the interval $[-2, 2]$ are used to determine the weights and biases in BI-blocks following equations 
	    \eqref{eq:bi-block:relu:weight:v2}, \eqref{eq:bi-block:relu:bias:v2} and \eqref{eq:bi-blocks:tanh:params:v2}.}
		\label{fig:bi-blocks:image}
	\end{center}
\end{figure}

Recall that the hyperbolic tangent function, $\tanh(x)$, is defined as the difference of two exponential functions divided by their sum. To simplify this expression, particularly near $x=0$, we apply the Taylor series expansion for the exponential terms. By employing a first-order approximation around $x=0$, we simplify $\tanh(x)$ as follows:
\ben
\tanh(x)= \frac{e^x-e^{-x}}{e^x+e^{-x}}
= \frac{1+x-(1-x)+o(x)}{1+x+(1-x)+o(x)}
= \frac{2x+o(x)}{2+o(x)}
\approx x,\quad x\rightarrow 0, 
\een
where $o(x)$ denotes higher-order terms that become negligible when $x$ becomes sufficiently small. Considering the asymptotic behavior of  $\tanh(x)$, which approaches $1$ as $x\rightarrow+\infty$ and $-1$ as $x\rightarrow-\infty$,  we introduce a piecewise linear function $s(x)$ that approximately captures the overall behavior of $\tanh(x)$:
\ben
s(x)=\left\{
\begin{array}{ll}
	-1,\quad x<-1,\\
	x,\qquad -1\leq x\leq 1,\\
	1,\qquad x>1.
\end{array}\right.\label{eq:f}
\een
This $s(x)$ provides an approximate simplification for the $\tanh(x)$ function, offering an accessible approximation for the construction of the BI-block with the Tanh activation function.
Indeed, the basis function $\varphi_j(x)$ can be reformulated as: 
\ben
\varphi_j(x)
=\frac{1}{2}s(\frac{2}{h_{j-1}}(x-x_{j-\frac{1}{2}}))-\frac{1}{2}s(\frac{2}{h_{j}}(x-x_{j+\frac{1}{2}})),
\een
with $x_{j-\frac{1}{2}}=x_j-\frac{1}{2}h_{j-1}$ and $x_{j+\frac{1}{2}}=x_j+\frac{1}{2}h_{j}$ being the midpoints of intervals $[x_{j-1},x_j]$ and $[x_j,x_{j+1}]$, respectively. 
Substituting $s(x)$ with $\tanh(x)$, we derive an approximating representation:
\ben
\varphi_j(x)
\approx\frac{1}{2}\tanh(\frac{2}{h_{j-1}}(x-x_{j-\frac{1}{2}}))-\frac{1}{2}\tanh(\frac{2}{h_{j}}(x-x_{j+\frac{1}{2}})).
\een
This formulation can be further expressed as a two-layer neural network:
\ben
\varphi_j(x)
\approx W_2\cdot (\tanh(W_1\cdot x+b_1)),
\een
where the weights and biases are defined by:
\ben
W_1=(\frac{2}{h_{j-1}},\frac{2}{h_{j}})^T,\qquad b_1=-(\frac{2}{h_{j-1}}x_{j-\frac{1}{2}},\frac{2}{h_{j}}x_{j+\frac{1}{2}})^T,\qquad W_2=(\frac{1}{2},-\frac{1}{2}).
\een
Similar to \eqref{eq:bi-block:relu:v2}, we rewrite the formulation as:
\be
\varphi_j(x)
\approx W_2\cdot\tanh(W_1^2\cdot(W_1^1\cdot x+b_1^1)),
\label{eq:bi-blocks:tanh:v2}
\ee
where
\be
W_1^1=(\sqrt{\frac{2}{h_{j-1}}},\sqrt{\frac{2}{h_{j}}})^T,\quad 
b_1^1=-(\sqrt{\frac{2}{h_{j-1}}}x_{j-\frac{1}{2}},\sqrt{\frac{2}{h_{j}}}x_{j+\frac{1}{2}})^T,\quad
W_1^2=\begin{bmatrix}
	\sqrt{\frac{2}{h_{j-1}}} & 0 \\
	0 & \sqrt{\frac{2}{h_{j}}}
\end{bmatrix}.
\label{eq:bi-blocks:tanh:params:v2}
\ee
Equation \eqref{eq:bi-blocks:tanh:v2} indicates that the basis function $\varphi_j(x)$ can be approximately represented as a three-layer network block as shown in Fig. \ref{fig:bi-blocks}(b), which we term the Basis-inspired Block with Tanh activation and denote by $\text{Block}_{\tanh}(x;x_j,h_{j-1},h_j)$. Similar to $\text{Block}_{\text{relu}}(x;x_j,h_{j-1},h_j)$, the weights and biases within block $\text{Block}_{\tanh}(x;x_j,h_{j-1},h_j)$ are also initialized using $x_j,h_{j-1}$ and $h_j$ according to equation  \eqref{eq:bi-blocks:tanh:params:v2}, which provides a good starting point for subsequent optimization. 
In Fig. \ref{fig:bi-blocks:image}, we depict the image of BI-blocks $\text{Block}_{relu}(x;x_j,h_{j-1},h_j)$ and $\text{Block}_{\tanh}(x;x_j,h_{j-1},h_j)$ where equally spaced nodes within interval $[-2,2]$ are used to compute the weights and biases according to Equations \eqref{eq:bi-block:relu:weight:v2}, \eqref{eq:bi-block:relu:bias:v2} and \eqref{eq:bi-blocks:tanh:params:v2}. 
\subsection{Basis-inspired DNNs}
With the BI-block in place, we are ready to introduce the novel Basis-inspired DNN (BI-DNN). For illustrative purposes, we first introduce the one-dimensional BI-DNN and then extend it to d-dimensional BI-DNN under the framework of the Kolmogorov Superposition Theorem.
\subsubsection{One-dimensional Basis-inspired DNN}
\begin{figure}[htbp]
	\begin{center}
		\includegraphics[width=0.85\textwidth]{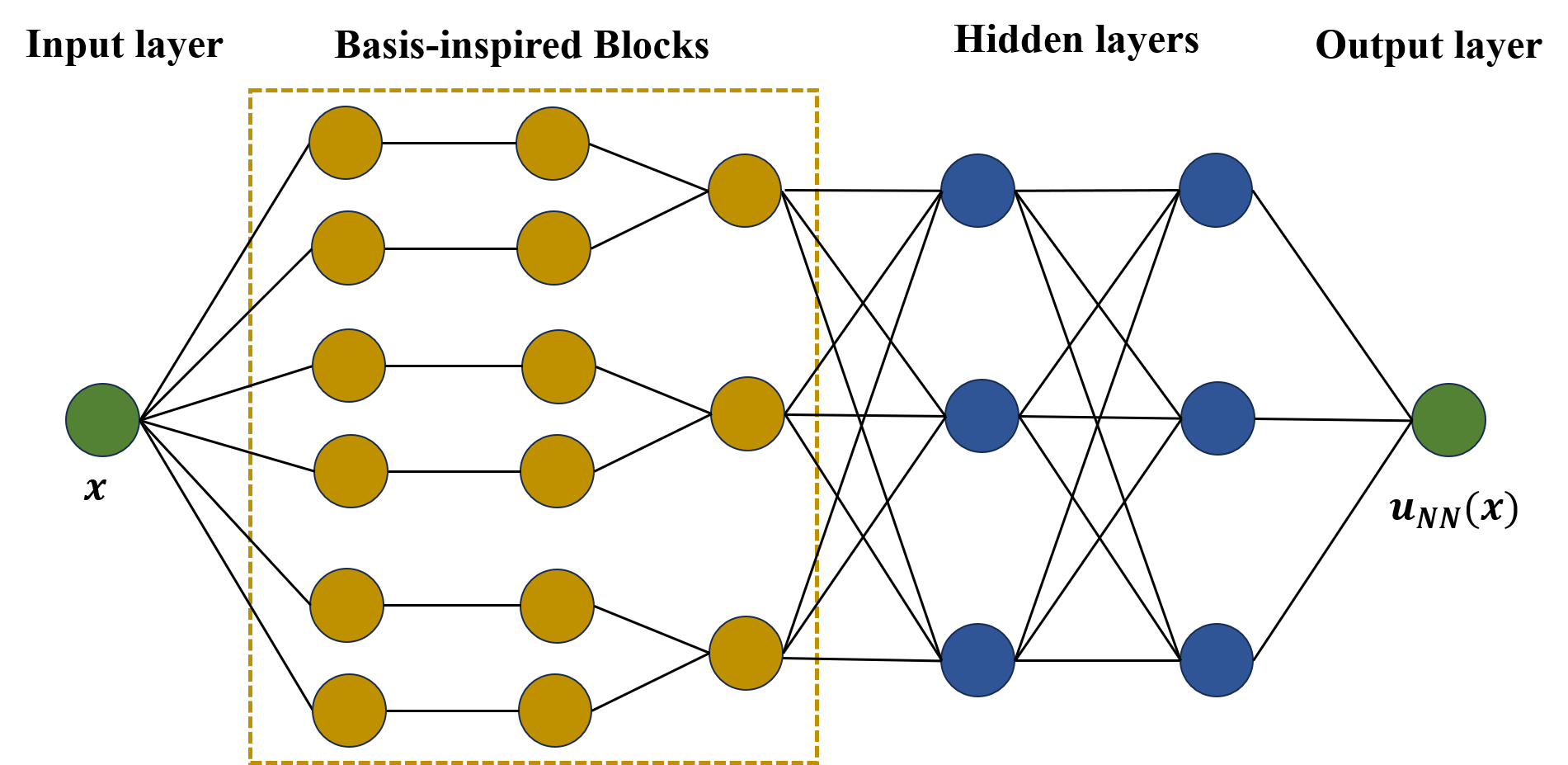}
	\end{center}
	\caption{Illustration of the architecture of a one-dimensional BI-DNN consisting of three BI-blocks with the Tanh activation followed by a three-layer subnetwork. The only difference to a conventional DNN is the  BI-blocks following the input layer.}
	\label{fig:bi-dnn:architecture:1d}
\end{figure}

Recall that the finite element function in conventional FEM is constructed as a linear combination of basis functions, that is,
\ben
u_{FEM}(x)=\sum_{j=1}^{n} u_j \varphi_j(x),\quad x\in I.
\een
It can be translated into the realm of neural networks by BI-blocks as:
\be
u_{NN}(x)=\sum_{j=1}^{n} u_j \text{Block}_{\text{relu}}(x;x_j,h_{j-1},h_j)
\approx \sum_{j=1}^{n} u_j \text{Block}_{\text{tanh}}(x;x_j,h_{j-1},h_j).
\label{eq:fem:function}
\ee
Motivated by expression \eqref{eq:fem:function}, we introduce the one-dimensional BI-DNNs as:
\be
b(x)=\mathbf{F}_{1:L}(\mathbf{B}(x)),
\label{eq:bi-dnn:1d}
\ee
where $\mathbf{F}_{1:L}$ is a fully connected subnetwork with $L$ layers expressed as:
\ben
\mathbf{F}_{1:L}=F_L \circ \sigma \circ F_{L-1} \circ \cdots \circ \sigma \circ F_{2} \circ \sigma  \circ F_1,
\een
and $\mathbf{B}$ is a stack of BI-blocks defined as:
\be
\mathbf{B}(x)=\left[B_1(x),B_2(x),\cdots,B_n(x)\right]^T.
\label{eq:bi-blocks:stacked}
\ee
For notational brevity, we utilize $B_j(x)$ in \eqref{eq:bi-blocks:stacked} to represent either $\text{Block}_{\text{relu}}(x;x_j,h_{j-1},h_j)$ or $\text{Block}_{\text{tanh}}(x;x_j,h_{j-1},h_j)$. It is easy to see that the expression \eqref{eq:fem:function} is a special BI-DNN where the fully connected subnetwork $\mathbf{F}_{1:L}$ is reduced to a single layer. An example of the architecture of a BI-DNN is depicted in Fig. \ref{fig:bi-dnn:architecture:1d}.

\par
To extend the one-dimensional BI-DNN to the multi-dimensional case, a direct method is to utilize tensor products to construct BI-blocks based on multi-dimensional basis functions. However, this approach encounters a significant challenge: the number of basis functions increases exponentially with dimension $d$, leading to the notorious curse of dimensionality (CoD). To circumvent this issue, 
we explore an alternative approach inspired by the Kolmogorov Superposition Theorem, which provides a framework for representing complex multivariate functions by univariate functions and offers a new insight into function approximation in high dimensions.
\subsubsection{Kolmogorov Superposition Theorem}
The Kolmogorov Superposition Theorem (KST), also known as Kolmogorov-Arnold representation theorem, was originally proposed by Kolmogorov \cite{kolmogorov1957representation} in 1957. The original version of the KST states that any function $f\in C(\left [ 0,1 \right ] ^{d})$ can be exactly represented by a finite composition of continuous univariate functions and the binary operator of addition.
After that, improved versions of the KST have emerged, focusing on reducing the number of univariate functions required and improving their smoothness. In particular, \cite{lorentz1962metric,lorentz1966approximation} propose a simplified version where the minimum feasible number of univariate functions are used. We state this version of KST as follows:  
\begin{thm}[Kolmogorov Superposition Theorem] 
	For any continuous function $f$ defined on $\left [ 0,1 \right ] ^{d}$, there exist irrational numbers $0< \lambda _{i} \le 1$ for $i=1,\,2,\,\cdots,\,d$, and strictly increasing $Lip \left ( \alpha  \right )$ inner functions $\phi_q$ (independent of $f$) with $\alpha = log_{10}2$ on $\left [ 0,1 \right ]$ for $q=0,\,1,\,\cdots,\,2d$, and with the presence of a continuous outer function $g(z),\,z \in \left [ 0,d \right ]$ such that the following identity holds:
	\begin{equation}
		f (x_{ 1},\cdots ,x _{d} ) = \sum_{q=0}^{2d} g\left ( \sum_{i=1}^{d} \lambda _{i} \phi _{q} \left ( x_{i} \right ) \right ).
		\label{eq:KST}
	\end{equation}
	\label{thm:KST}
\end{thm}
For clarify, the class $\mathrm{Lip}(\alpha),0 < \alpha \leq 1$ consists of all functions that satisfy $|f(x) - f(y)| \leq M|x - y|^\alpha,\; \forall\,x,\,y \in [a, b]$ for some $M$.
Equation \eqref{eq:KST} reveals that the task of learning a multivariate function is simplified to learning $2d+2$ univariate functions. Unfortunately, the construction of these univariate functions involves an infinite process 
\cite{Sprecher_Draghici_2002}, 
which results in the unachievability of exact representations through the KST in practice. Hence, researchers have increasingly focused on exploring approximate versions of KST where univariate functions involved are approximated with simpler and smoother functions, such as splines function \cite{lai2023kolmogorov} and neural networks \cite{Kurkova_1992, Schmidt-Hieber_2020}. In the spirit of approximate versions of KST, we construct our multi-dimensional BI-DNN by utilizing specially designed subnetworks to approximate univariate functions and assemble them according to the KST framework.

\subsubsection{Multi-dimensional Basis-inspired DNNs}
\begin{figure}[H]
	\centering
	\includegraphics[width=0.8\textwidth]{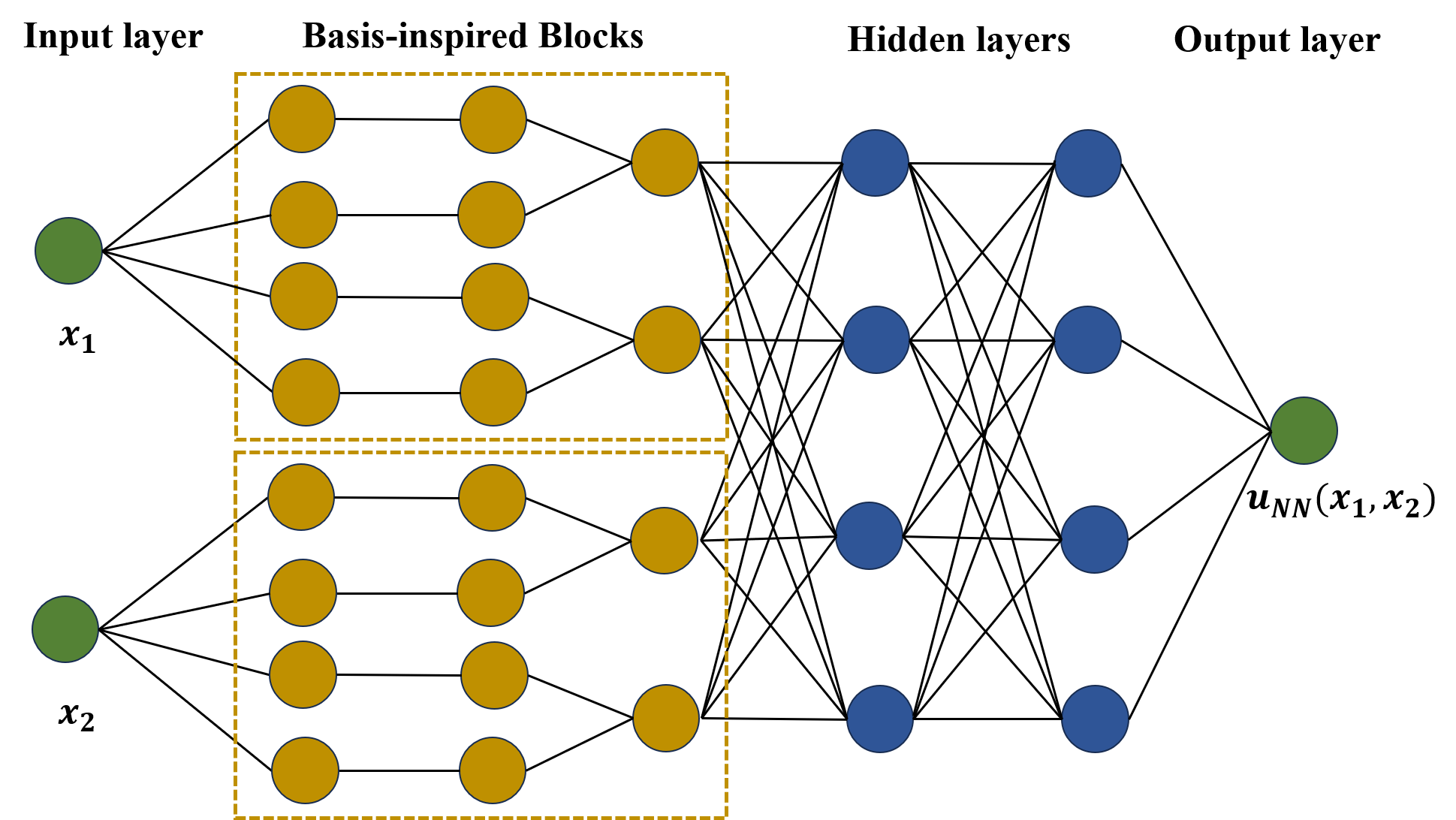}
	\caption{Illustration of the architecture of a two-dimensional BI-DNN consisting of two BI-blocks with Tanh activation function in each dimension.} 
	\label{fig:bi-dnn:architecture:2d}
\end{figure}

For clarity, we adopt the vectorized form of representation previously formulated in \cite{Zhang_Li_Chang_Liu_Huang_Xiang_2022}. For any $x\in [0,1]$ and $\boldsymbol{z}=(z_0,z_1,\cdot\cdot\cdot,z_{2d})\in[0,d]^{2d+1}$, we define a univariate vector-valued function  $\boldsymbol{\Phi}(x)$ and a multivariate vector-valued function $\mathbf{G}(\boldsymbol{z})$ as follows:
\begin{equation*}
	\boldsymbol{\Phi}(x)=(\phi_0(x),\phi_1(x),\cdot\cdot\cdot,\phi_{2d}(x)),\quad \mathbf{G}(\boldsymbol{z})=(g(z_0),g(z_1),\cdot\cdot\cdot,g(z_{2d}))^T.
\end{equation*}
Equation \eqref{eq:KST} can be rewritten as: 
\begin{equation}
	f(x_1,\cdot\cdot\cdot,x_d)=V\cdot \mathbf{G}(\Lambda \cdot [\boldsymbol{\Phi}(x_1),\boldsymbol{\Phi}(x_2),\cdot\cdot\cdot,\boldsymbol{\Phi}(x_d)]^T),
	\label{eq:KST_vec}
\end{equation}
where $\Lambda=[\lambda_1 \mathbf{I}_{2d+1},\lambda_2 \mathbf{I}_{2d+1},\cdot\cdot\cdot, \lambda_d \mathbf{I}_{2d+1}]$ with $\mathbf{I}_{2d+1}$ being a $(2d+1)\times(2d+1)$ identity matrix, and $V=[1,1,\cdot\cdot\cdot,1]$ being an all-one row vector of length $2d+1$. Notice that $\Lambda$ and $V$ represent linear transformations, which, within the context of the neural network, can be viewed as affine transforms $F_1$ and $F_2$, respectively. To approximately reformulate expression \eqref{eq:KST_vec} into the format of neural network, we employ one-dimensional BI-DNNs defined in \eqref{eq:bi-dnn:1d} to approximate $\boldsymbol{\Phi}(x_j),j=1,2\cdots,d$ and a fully connected subnetwork $\mathbf{F}_{1:L_G}$ to approximate $\mathbf{G}$.
Consequently, the continuous multivariate function can be approximated in the neural network architecture as follows:
\be
f(x_1,\cdot\cdot\cdot,x_d) \approx F_2 \circ  \mathbf{F}_{1:L_G}\circ F_1([\mathbf{F}_{1:L_1}(\mathbf{B_1}(x_1)),\mathbf{F}_{1:L_2}(\mathbf{B_2}(x_2)),\cdots,\mathbf{F}_{1:L_d}(\mathbf{B_d}(x_d))]^T),
\label{eq:KST_vec:approx}
\ee
where $\mathbf{B}_j$ is a collection of BI-blocks as represented in \eqref{eq:bi-blocks:stacked} and  $\mathbf{F}_{1:L_1},\cdots,\mathbf{F}_{1:L_d}$ denote distinct fully connected subnetworks. Motivated by \eqref{eq:KST_vec:approx}, we define the multi-dimensional BI-DNN as: 
\be
b(x_1,x_2,\cdots,x_d)=\mathbf{F}_{1:L}([\mathbf{B}_1(x_1),\mathbf{B}_2(x_2),\cdots,\mathbf{B}_d(x_d)]).
\label{eq:bi-dnn:dd}
\ee
We remark that the network architecture defined in \eqref{eq:bi-dnn:dd} is a generalization of the one outlined in  \eqref{eq:KST_vec:approx}. In this generalization, sparse connections among  $F_{1:L_1},\cdots,F_{1:L_d}$
are generalized to fully connected layers, and activation functions are incorporated into previously inactivated hidden layers. This generalization leads to a more simplified and powerful model to approximate complex multivariate functions.
An example of the architecture of a BI-DNN with $d=2$ is shown in Fig. \ref{fig:bi-dnn:architecture:2d}.
As we conclude this section, it is worth noting that the weights and biases of BI-blocks are closely related to the nodal positions, a relationship particularly evident when the ReLU activation function is employed (as shown in equations \eqref{eq:bi-block:relu:weight:v2} and \eqref{eq:bi-block:relu:bias:v2}). As the training of the  BI-DNN involves optimizing weights and biases in BI-blocks, it implies initial nodal coordinates are indirectly optimized as well. This approach is similar in spirit to the r-adaptivity in FEM, where nodes of the mesh elements are dynamically moved to better capture the solution of the PDE of interest. In addition, h-adaptivity is also a favorable adaptive technique in FEM that achieves high accuracy with an economical use of degrees of freedom by locally refining the mesh and introducing additional basis functions where needed. Motivated by h-adaptivity, we further study the idea of incrementally adding BI-blocks focusing on local regions with high errors.

\section{Adaptive Basis-inspired DNN}
In this section, we present the Adaptive Basis-inspired Deep Neural Network (ABI-DNN), an adaptive framework designed to incrementally generate an appropriate BI-DNN for a specific learning task.  

Note that standard AFEMs employ iterative cycles of `solve, estimate, mark, refine' to ensure the pre-specified tolerance is satisfied. Following a similar strategy, we start with a small size BI-DNN, solve the optimization problem in \eqref{eq:loss}, and estimate the error distribution by the local error indicator
\be
\eta_p=\left|\mathcal{L}u_{NN}(x_p)-f(x_p))\right|, \quad \text{for}\;\;\; x_p\in S,
\label{eq:indicator}
\ee
where $S$ is a set of training samples randomly selected in advance. To identify spatial regions with large errors, we first mark points with large errors by the maximum marking strategy:
\be\label{eq:mark:maximum}
\hat{S}=\left\{x_p\in S:\eta_p>\gamma\max_{x_p\in S}\{\eta_p\},\qquad \text{for}\ \gamma\in (0,1)\right\}.
\ee
Alternative marking strategies, such as the bulk marking strategy \cite{ZhimingChen2010Selected}, are also available.
\par
After identifying points with large errors, we aim to enhance the network architecture by integrating new BI-blocks to better resolve local regions with significant errors. As shown in Equations \eqref{eq:bi-block:relu:v2} and \eqref{eq:bi-blocks:tanh:v2}, the initialization of a new BI-block is determined by $x_j,h_{j-1}$ and $h_j$. These parameters determine where to place the newly added BI-block, and thus, they should be carefully chosen to guide the new BI-block to focus on a specific local region. An effective approach to address the issue is to cluster the marked points in $\hat{S}$ and then locate BI-blocks at the centroid of each cluster. Various clustering techniques are available for this purpose. Given the challenge of predetermining the number of clusters, density-based clustering, such as DBSCAN \cite{2012DBSCAN}, is a good choice. Once the clusters are identified, the coordinate of the centroid and radius of each cluster can be calculated, and then BI-blocks are added accordingly. To be specific, if the centroid of a cluster is $(c_1,c_2,\cdots,c_d)$ and the radius is $r$, then BI-blocks:
\be
\text{Block}_{\sigma}(x_1;c_1,r_s,r_s),\quad \text{Block}_{\sigma}(x_2;c_2,r_s,r_s),\quad \cdots,\quad \text{Block}_{\sigma}(x_d;c_d,r_s,r_s),
\label{eq:cluster:new:bi-block}
\ee
with $r_s = s\cdot r$ and $s$ being a scaled constant to enhance the flexibility, are added to
\ben \mathbf{B}_1(x_1),\quad
\mathbf{B}_2(x_2),\quad
\cdots,\quad
\mathbf{B}_d(x_d), 
\een 
in Equation \eqref{eq:bi-dnn:dd}, respectively. Consequently, the number of clusters identified will determine the number of blocks added to each $\mathbf{B}_j(x_j),j=1,2,\cdots,d$. To avoid potential information bottlenecks due to drastic changes in layer width, it is beneficial to maintain a gradual transition in layer widths between BI-blocks and subsequent hidden layers. A heuristic strategy is periodically adding new neurons to the subsequent hidden layers as the number of BI-blocks increases. After the network architecture is properly enhanced, we proceed to solve the optimization problem \eqref{eq:loss} using the updated network architecture. We repeat the above procedure until the total error indicator $\eta=\left(\frac{1}{N}\sum_{x_p\in\mathcal{S}}\eta_p^2\right)^{1/2}$ meets the predefined tolerance criterion or the predetermined maximum number of enhancement iterations is reached. The ABI-DNN method is summarized in Algorithm \ref{algo:ABI-DNN}.
\begin{algorithm}[htbp]
	\caption{ABI-DNN Method}
	\label{algo:ABI-DNN}
	\KwIn{A tolerance $\eta_{tol}$, the maximum number of adaptive iterations $J$, an initial small-sized BI-DNN, a training sample set $S$}
	\KwOut{The numerical solution $u_{NN}$ and the automatically generated ABI-DNN}
	\BlankLine
	Solve the optimization problem \eqref{eq:loss};\\
	Estimate the local indicator $\eta_p$ 
	and compute the total error indicator $\eta$;\\
	\While{$\eta>\epsilon$ \text{and the iteration count is less than} $J$}{
		Mark points in $S$ by \eqref{eq:mark:maximum} to generate $\hat{\mathcal{S}}$;\\ Cluster points in $\hat{\mathcal{S}}$ into distinct groups $\mathcal{C}$;\\
		Refine the current BI-DNN according to clusters $\mathcal{C}$;\\
		Solve the optimization problem \eqref{eq:loss} with the updated BI-DNN;\\
		Re-estimate the local indicator $\eta_p$ and recompute the total esitimator $\eta$;\\
	}
\end{algorithm}

Finally, we detail the initialization strategies utilized in the ABI-DNN method.
It is well known that a good initialization is crucial to obtaining an ideal approximate solution, especially in the context of nonlinear optimizations, which are known to yield multiple solutions. Fortunately, the ABI-DNN framework naturally comes with a natural process for attaining a good initialization. Specifically, at the beginning of the ABI-DNN method, 
we first subdivide the computational domain into a coarse uniform mesh, then BI-blocks in the initial BI-DNN are initialized utilizing the coordinates of nodes on the mesh according to Equation \eqref{eq:bi-block:relu:v2} or \eqref{eq:bi-blocks:tanh:v2}. Weights in subsequent hidden layers are initialized using the Xavier initialization method. During the adaptive iteration, the approximation achieved by the preceding BI-DNN is already a good approximation to the current BI-DNN. Therefore, the parameters of the existing neurons inherit the current approximation as their initial values, and the initialization of each newly added BI-block is guided by the centroid and the radius of its associated cluster, as shown in Equation  \eqref{eq:cluster:new:bi-block}. To initialize the weights and biases of newly added neurons that are not encapsulated within any BI-block, a straightforward method is to set them to zero. This approach ensures that the initial approximation is the current approximation, providing a good starting point for subsequent training.

\section{Numerical Experiments}
In this section, we present a series of numerical experiments to demonstrate the superior performance of the proposed methods. In particular, Subsection 5.2 explores function fitting for both singular and smooth cases. Subsection 5.3 investigates Poisson equations, focusing on solutions with peaks and domains with re-entrant corners. Subsection 5.4 studies Burger's equation. 

\subsection{Experiment Setup}
In all experiments, we use standard PINN as the baseline for comparison. For simplicity, the hidden layers of the underlying neural network are chosen as fully connected layers with an equal number of neurons. For the sake of clarity, we use PINN($w=m$) or DNN($w=m$) to refer to the underlying network that has $m$ neurons in each hidden layer. In addition, we introduce BI-DNN($b=[m_1,\cdots,m_d]$) and ABI-DNN($b=[m_1,\cdots,m_d]$) to denote the corresponding neural networks consist of $m_1,\cdots,m_d$ BI-blocks in each dimension, respectively. In the case of $d=1$, these notations are simplified to BI-DNN($b=m$) and ABI-DNN($b=m$).
\par
As illustrated in the Algorithm \ref{algo:ABI-DNN}, the ABI-DNN method starts at a small-size BI-DNN with initialization described in Section 4. During the adaptive learning, the adaptive model is optimized $10000$ epochs for one-dimensional problems and $15000$ epochs for two-dimensional problems after each adaptive enhancement. 
To spot local regions requiring enhancement, local indicator \eqref{eq:indicator} is evaluated on the predefined training set, and maximum marking strategy \eqref{eq:mark:maximum} with $\gamma=0.5$ is used to mark points with large errors. 
To circumvent the challenge of presetting cluster numbers, we employ the widely recognized density-based clustering algorithm DBSCAN to cluster the marked points. We remark that DBSCAN requires the specification of two critical parameters: the maximum distance, denoted as $\epsilon$, between two samples for one to be considered in the neighborhood of the other, which we set to $\epsilon=0.1$; and the minimum number of neighbors, referred to as $MinPts$, within the $\epsilon$-radius to form a dense region, which we set to $MinPts=1$ to prevent any marked point from being labeled as noise. In addition, the $L_\infty$ norm is adopted to measure the distance between pairs of points. 
Once marked points are clustered, new BI-blocks as defined in  \eqref{eq:cluster:new:bi-block} with $s=2$ are incorporated into the current BI-DNN, and the updated BI-DNN will be continuously trained.
To avoid excessive computational costs, the adaptive iteration process will stop once the maximum number of iteration steps $J$ is reached or the prescribed tolerance $\eta_{tol}$ is obtained. For all subsequent experiments, the maximum number of iterations is set to $J=10$, and the tolerance $\eta_{tol}$ will be explicitly specified in each experiment.
\par
In all experiments, the minimization problems are iteratively solved by the Adam optimizer\cite{kingma2014adam}. The learning rate starts with an initial value of $\tau = 5 \times 10^{-3}$ and decays every $2500$ steps with a base of $0.9$. The penalty parameter in loss function \eqref{eq:loss} is set to $\beta=1000$ in all experiments.
To assess the accuracy of the learned solution, we define the relative $L_2$ error:
\ben
\text{Error} = \frac{\sqrt{\sum_{n=1}^N |u_{NN}(x^{(n)},\theta)-u_{*}(x^{(n)})|^2}}{\sqrt{\sum_{n=1}^N |u_{*}(x^{(n)})|^2}},
\een
where $u_*$ denote the target function, $u_{NN}$ denote the approximating function learned by PINN, BI-DNN or ABI-DNN, and $\{x^{(n)}\}_{n=1}^N$ are $N$ uniformly distributed testing data points with $N=500$ in one-dimensional case and $N=40000$ in the two-dimensional case unless otherwise stated.
\subsection{Function Fitting}
To investigate the effectiveness of the proposed ABI-DNN method, we first consider two function-fitting problems. The first test problem is a singular target function:
\begin{equation}
	u_*(x) = 
	\begin{cases} 
		\;	25x^2,&\quad x\in[0,\;0.2),\\
		\;	25(0.4-x)^2,&\quad x\in[0.2,\;0.4), \\
		\;	0,&\quad x\in[0.4,\;1],
	\end{cases}
	\label{pb:1d:singular}
\end{equation}
which is characterized by the presence of a cusp at $x=0.2$, indicating a point of non-differentiability that poses a challenge for numerical approximations (see Fig. \ref{fig:fitting1d:exact:singular}(a) for illustration). The second test problem is a smooth high-frequency target function (see Fig. \ref{fig:fitting1d:exact:smooth}(a) for illustration): 
\be
u_*(x)=\sum_{i=0}^5 \sin(2^i\pi x), \quad x\in [0,1].
\label{pb:1d:smooth}
\ee
For the sake of simplicity, we only consider BI-DNNs and ABI-DNNs that include no additional hidden layers beyond those containing BI-blocks. Specifically, the output layer immediately follows the BI-blocks without any intervening hidden layers. Thus, once the number of BI-blocks is chosen and the activation function is determined, the architecture of the BI-DNN is established. Note that the BI-DNN, configured in this manner, consists of four layers in total. When comparisons with DNNs are necessary, four-layer DNNs with an equal number of neurons in each hidden layer are employed to ensure a fair comparison. The loss function is defined as the discrete mean square error 
\begin{equation*}
	L(\theta) = \frac{1}{N_{r}} \sum_{n=1}^{N_{r}}|u_{NN}(x^{(n)};\theta)-u_*(x^{(n)})|^{2},
\end{equation*}
where $N_{r} = 2000$ points are randomly sampled  from the interval $[0,1]$ to form the training dataset. 
\subsubsection{The effect of r-adaptivity}
\begin{figure}[htbp]
	\begin{center}
		\begin{tabular}{cccc}
			\centering
			\includegraphics[width=0.4\textwidth]{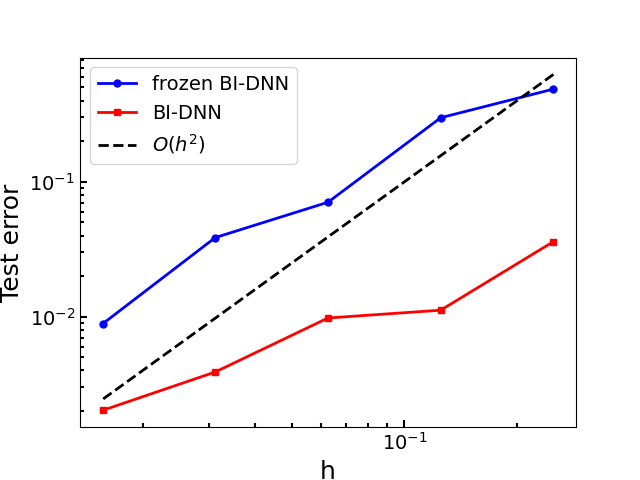} &
			\includegraphics[width=0.4\textwidth]{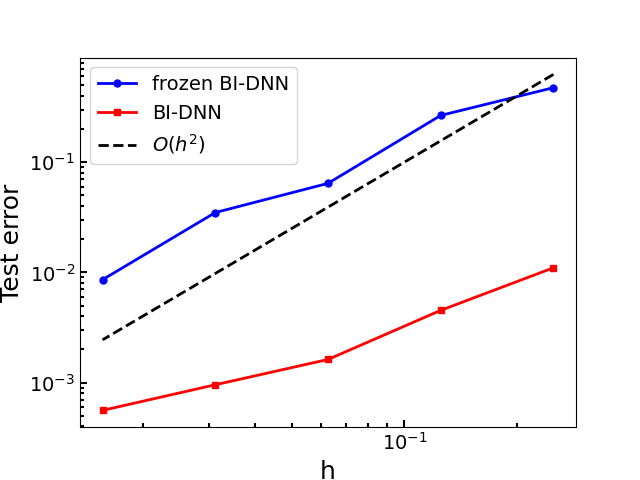} \\
			(a) ReLU activation, Problem \eqref{pb:1d:singular} & (b) Tanh activation, Problem \eqref{pb:1d:singular} \\
			\includegraphics[width=0.4\textwidth]{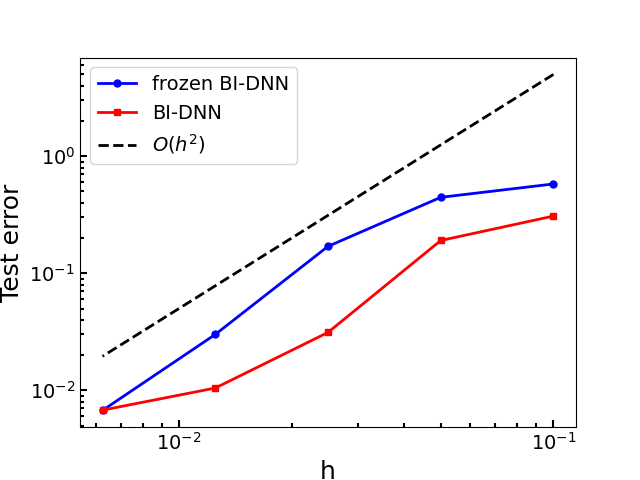} &
			\includegraphics[width=0.4\textwidth]{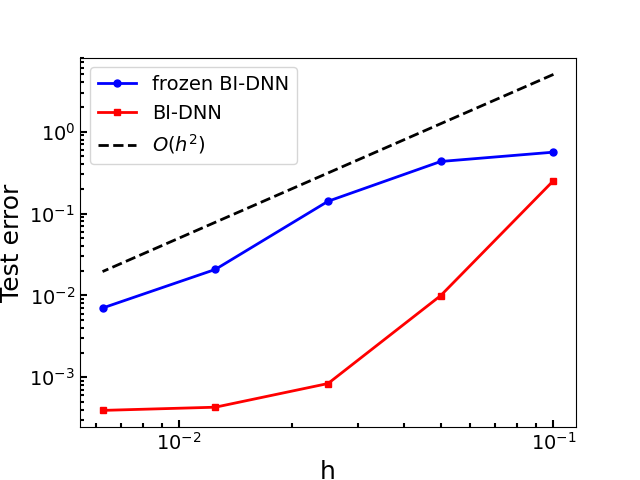}\\
			(c) ReLU activation, Problem  \eqref{pb:1d:smooth} & (d) Tanh activation, Problem \eqref{pb:1d:smooth}
		\end{tabular}
	\end{center}
	\caption{Comparison of relative errors between frozen BI-DNNs and BI-DNNs on Problem \eqref{pb:1d:singular} (Top row) and Problem \eqref{pb:1d:smooth} (Bottom row).} 
	\label{fig:fitting1d:frozen}
\end{figure}

As detailed in Section 3.1, weights and biases of BI-blocks are closely related to the nodal positions. Consequently, the training process of BI-DNN involves indirect adjustments of nodal coordinates, which is similar in spirit to the r-adaptivity in FEM. For ease of reference, we use `frozen BI-DNN' to denote the BI-DNN  where the weights and biases within its BI-blocks are fixed, which implies that the nodal positions used to initialize its BI-blocks remain constant throughout the training process of the frozen BI-DNN. To investigate the convergence rate of the frozen BI-DNN and BI-DNN(with inherent r-adaptivity), the computational interval $[0,1]$ is uniformly partitioned into $m$ sub-intervals with $m=4, 8, 16, 32, 64$. The nodal coordinates resulting from this partition are used to initialize the weights and biases of the BI-blocks, which means BI-DNNs and frozen BI-DNNs consisting of $5, 9, 17, 33,$ and $65$ BI-blocks are utilized. 
After 50000 epochs of training, the results of the $L^2$ relative error are plotted in Fig. \ref{fig:fitting1d:frozen}. As expected, frozen BI-DNNs have the same convergence rates as linear FEMs. While BI-DNNs seem to exhibit slower convergence rates compared to frozen BI-DNNs, they achieve significantly higher accuracy with the same number of nodes. In addition, BI-DNNs with Tanh activation function demonstrate superior efficiency over those with ReLU activation function. Therefore, BI-DNNs with Tanh activation function are adopted in subsequent experiments.
\subsubsection{The effect of BI-DNN}
\begin{figure}[htbp]
	\begin{center}
		\begin{tabular}{cc}
			\centering
			\includegraphics[width=0.45\textwidth]{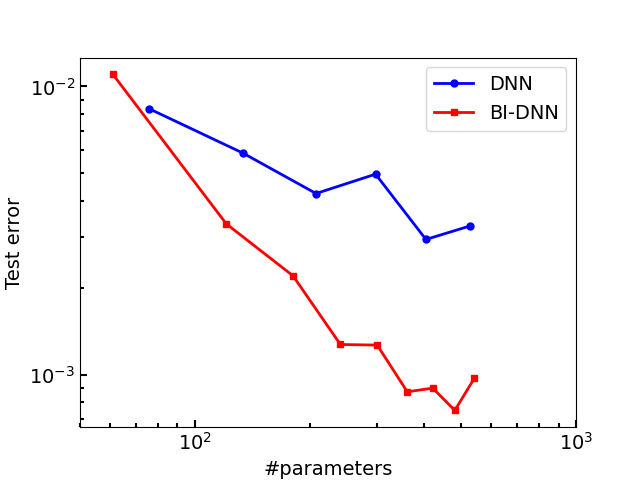} &
			\includegraphics[width=0.45\textwidth]{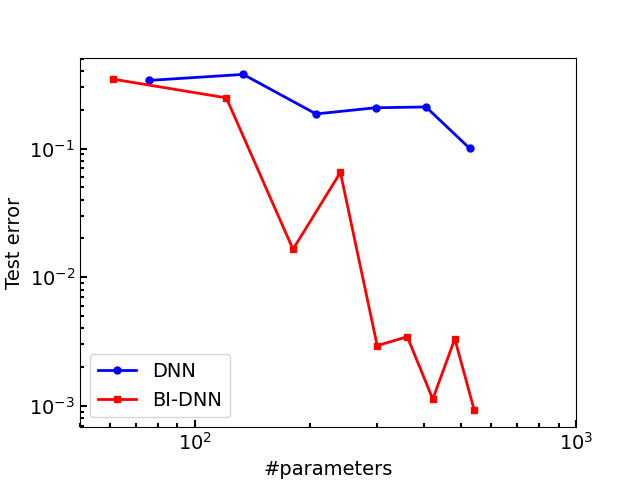} \\
			(a) & (b)
		\end{tabular}
	\end{center}
	\vspace{-0.2in}
	\caption{The change of relative errors of DNN and BI-DNN along with the number of trainable parameters on Problem \eqref{pb:1d:singular} (a) and Problem \eqref{pb:1d:smooth} (b), respectively.}
	\label{fig:fitting1d:params}
\end{figure}
To further evaluate the performance of BI-DNNs, a comparative analysis is conducted with conventional DNNs. In this experiment, BI-DNNs with $5,10,15,20,25,30,35,40,45$ BI-blocks are trained to study their performance as the number of trainable parameters increases. Comparatively, DNNs with 
$w=5,7,9,11,13,15$ neurons in each hidden layer were trained. After 50,000 epochs of training, the change of relative errors along with the number of parameters is illustrated in Fig. \ref{fig:fitting1d:params}. As expected, both the DNN and the BI-DNN exhibit decreasing error trends as the number of parameters increases, yet the BI-DNN demonstrates a faster rate of decrease. Furthermore, BI-DNNs consistently achieve significantly lower relative errors compared to DNNs with a similar number of trainable parameters. In particular, on the problem \eqref{pb:1d:smooth}, a BI-DNN with $541$ trainable parameters can achieve a relative error of $9.19E-04$, which is two orders of magnitude lower than the counterpart DNN with $526$ parameters, reaching an error of $1.00E-01$. These observations demonstrate that the BI-DNN can outperform the conventional DNN, particularly when addressing problems with localized features.
\subsubsection{The effect of ABI-DNN}
In this subsection, We further compare the performance of ABI-DNN with network models of fixed structure on problems \eqref{pb:1d:singular} and  \eqref{pb:1d:smooth}.

\begin{figure}[htbp]
	\begin{center}
		\begin{tabular}{ccc}
			\centering
			\includegraphics[width=0.3\textwidth]{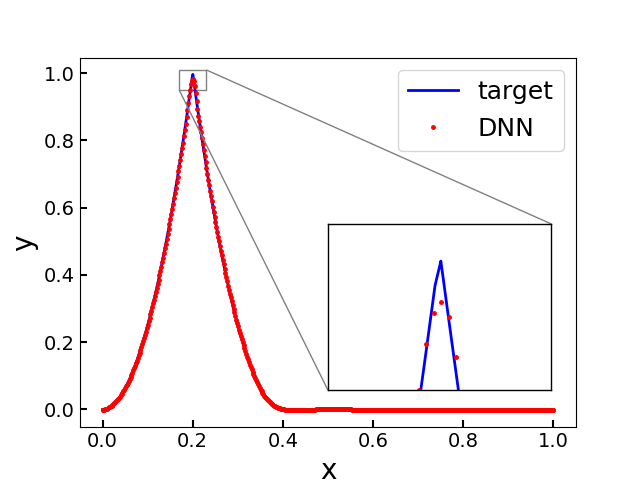} &
			\includegraphics[width=0.3\textwidth]{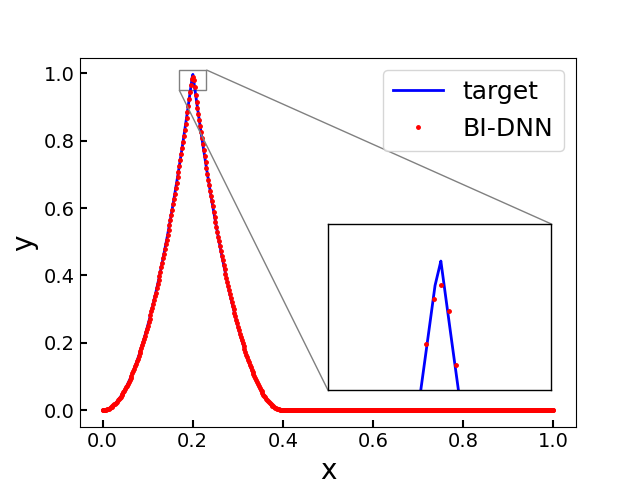} &
			\includegraphics[width=0.3\textwidth]{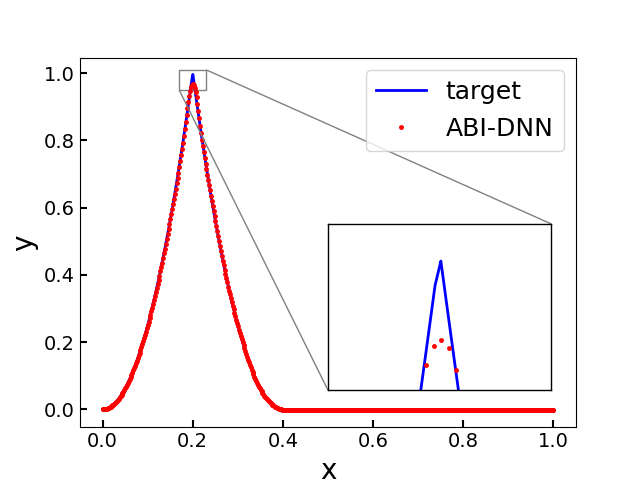}\\
			(a) DNN & (b) BI-DNN & (c) ABI-DNN, 1st iteration\\
			\includegraphics[width=0.3\textwidth]{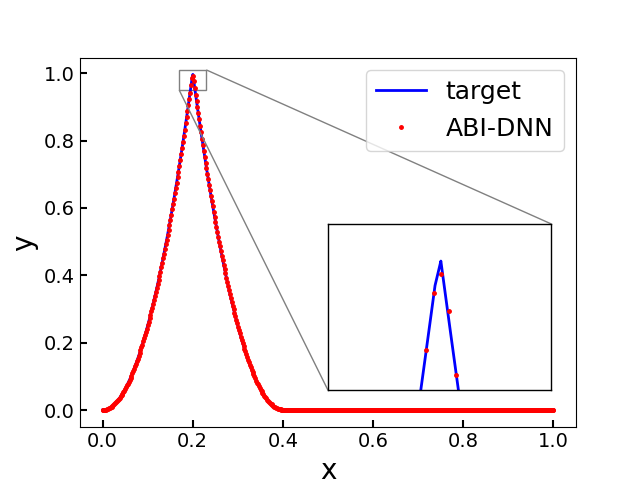} &
			\includegraphics[width=0.3\textwidth]{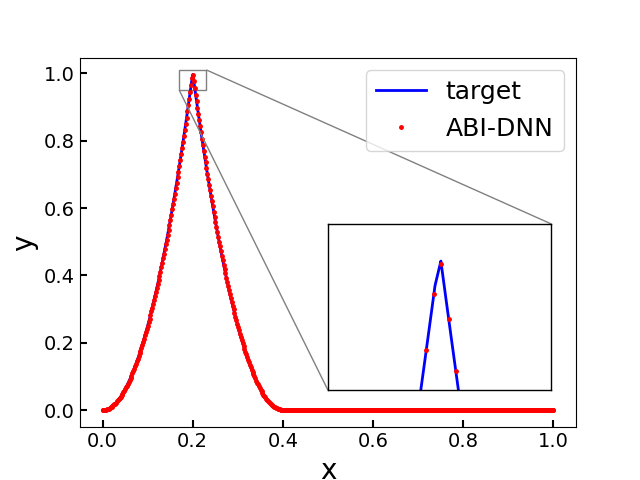} &
			\includegraphics[width=0.3\textwidth]{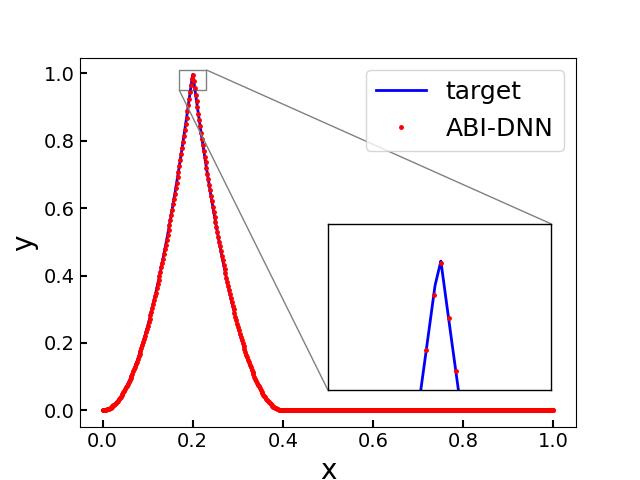} \\
			(d) ABI-DNN, 2nd iteration & (e) ABI-DNN, 4th iteration & (f) ABI-DNN, 5th iteration\\
		\end{tabular}
	\end{center}
	\caption{Illustration of the target function and approximating functions on problem \eqref{pb:1d:singular}. Blue curves in (a-f) represent the target function, and the red dots represent approximating functions given by: (a) DNN(w=9), (b) BI-DNN(b=16), (c-f) ABI-DNN after 1st, 2nd, 4th and 5th adaptive iterations, respectively. The insets of each subfigure are zoomed views of approximating solutions within the region $[0.17, 0.23]\times [0.95, 1.01]$.}
	\label{fig:fitting1d:exact:singular}
\end{figure}

\begin{figure}[htbp]
	\begin{center}
		\begin{tabular}{ccc}
			\centering
			\includegraphics[width=0.3\textwidth]{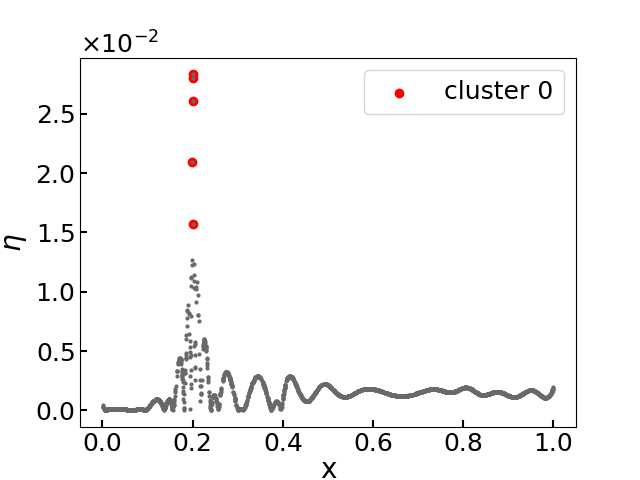}&
			\includegraphics[width=0.3\textwidth]{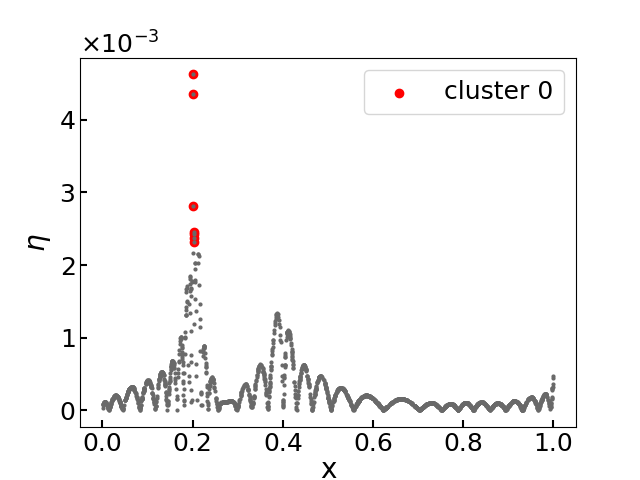}&
			\includegraphics[width=0.3\textwidth]{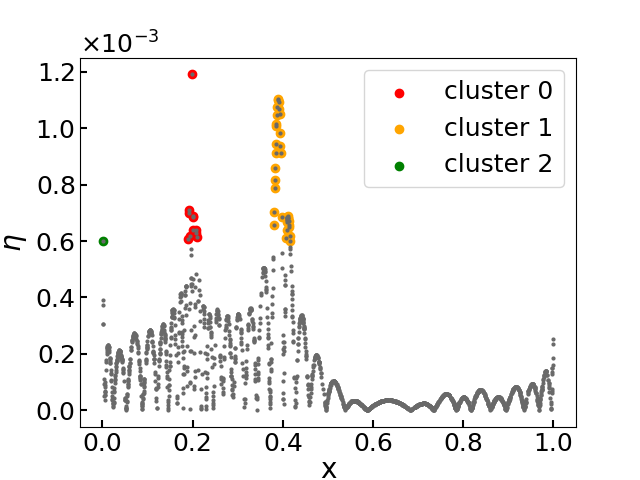}\\
			(a) after 1st iteration  & (b)  after 2nd iteration & (c)  after 4th iteration
		\end{tabular}
	\end{center}
	\caption{Illustration of clustering results of ABI-DNN after 1st, 2nd, and 4th adaptive iterations on Problem \eqref{pb:1d:singular}.}
	\label{fig:fitting1d:cluster:singular}
\end{figure}

\begin{figure}[htbp]
	\begin{center}
		\begin{tabular}{cc}
			\centering
			\includegraphics[width=0.45\textwidth]{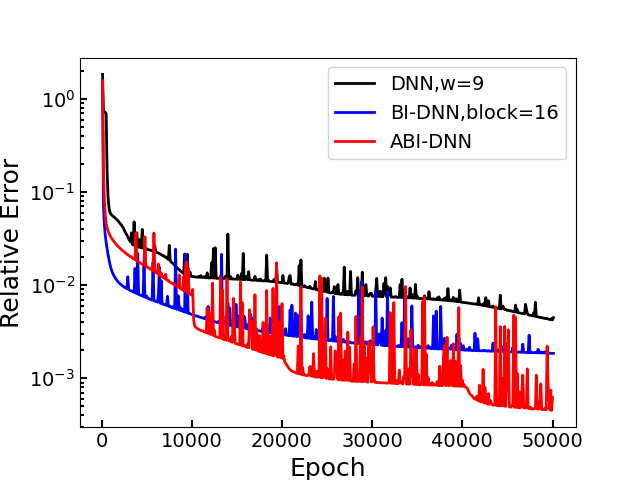}&
			\includegraphics[width=0.45\textwidth]{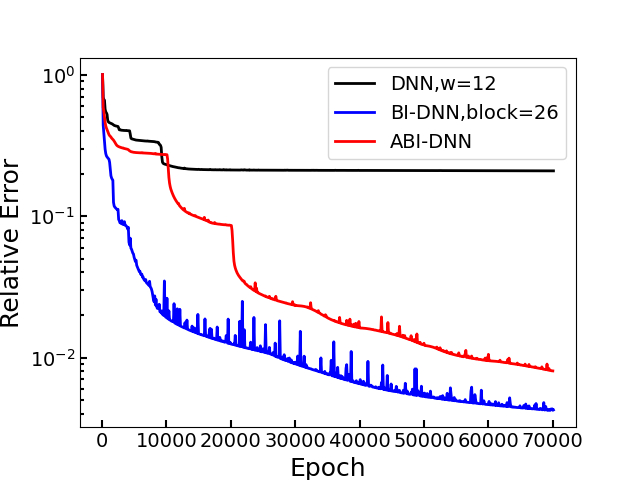}\\
			(a) Problem \eqref{pb:1d:singular} & (b) Problem \eqref{pb:1d:smooth}
		\end{tabular}
	\end{center}
	\caption{Convergence processes of DNN, BI-DNN and ABI-DNN on Problem \eqref{pb:1d:singular} (a) and Problem \eqref{pb:1d:smooth} (b).}
	\label{fig:fitting1d:convergence}
\end{figure}

For the singular target function \eqref{pb:1d:singular}, we start with a small size BI-DNN containing $10$ BI-blocks and set the tolerance $\eta_{tol}=0.0002$ in this example. After $5$ times of adaptive iteration, the ABI-DNN stops at a BI-DNN with $16$ BI-blocks, achieving a relative error of $4.51E-04$. The relative errors obtained in each iteration and the corresponding parameter counts used are reported in Table \ref{tab:fitting1d:ABI-DNN:singular}. For comparison, 
the BI-DNN($b=16$) and DNN($w=9$) are both trained for $50000$ epochs to match the number of training epochs of the ABI-DNN. The results are presented in the last two rows of Table  \ref{tab:fitting1d:ABI-DNN:singular}. It is clear that the final adaptive model ABI-DNN($b=16$) yields an error that is one order of magnitude lower than the errors obtained by
the BI-DNN($b=16$) and DNN($w=9$), both of which have a comparable number of parameters.  
Fig. \ref{fig:fitting1d:exact:singular} presents a comparison of the approximations provided by these three models. Upon examining the zoomed-in views in each subfigure, it is evident that DNN, BI-DNN and the ABI-DNN after first iteration all struggle to accurately capture the sharp feature at $x=1/5$. However, as adaptive enhancements progress, the ABI-DNN gives increasingly accurate representations of this sharp feature.
A possible reason for this improvement is the ABI-DNN's concentrated focus on addressing the singularity. As depicted in Fig. \ref{fig:fitting1d:cluster:singular}, points with large errors are marked and clustered, which guides the progressive addition of new BI-blocks around the singularity and leads to a more precise approximation.
The convergence behavior of these three models is illustrated in Fig. \ref{fig:fitting1d:convergence}(a). Notably, the ABI-DNN shows a significant reduction in error after each adaptive enhancement, ultimately achieving a lower final error than fixed models. This highlights the superiority of adaptive architecture enhancement for accurately characterizing singular features. 
\begin{table*}[htbp]
	\centering
	\caption{Numerical results of ABI-DNN on Problem \eqref{pb:1d:singular}. For comparison, the results of BI-DNN and DNN with a similar number of parameters are shown in the last two rows.
	}
	\begin{tabular}{c|ccccc}
		\toprule
		Model & Network structure & $\sharp$  Parameters & Testing error \\
		\midrule
		ABI-DNN(b=10) & 1-20-20-10-1 & 121 & 7.85E-03  \\
		\midrule
		ABI-DNN(b=11) & 1-22-22-11-1 & 133 & 1.63E-03  \\
		\midrule
		ABI-DNN(b=12) & 1-24-24-12-1 & 145 & 9.08E-04  \\
		\midrule
		ABI-DNN(b=13) & 1-26-26-13-1 & 157 & 8.19E-04  \\
		\midrule
		\textbf{ABI-DNN(b=16)} & \textbf{1-32-32-16-1} & \textbf{193} & {\textbf{4.51E-04}} \\
		\midrule
		BI-DNN(b=16)  & 1-32-32-16-1 & 193 & 1.85E-03  \\
		\midrule
		DNN(w=9)          & 1-9-9-9-1    & 208 & 4.25E-03  \\
		\bottomrule
	\end{tabular}
	\label{tab:fitting1d:ABI-DNN:singular}
\end{table*}

\par
\begin{figure}[htbp]
	\begin{center}
		\begin{tabular}{ccc}
			\centering
			\includegraphics[width=0.3\textwidth]{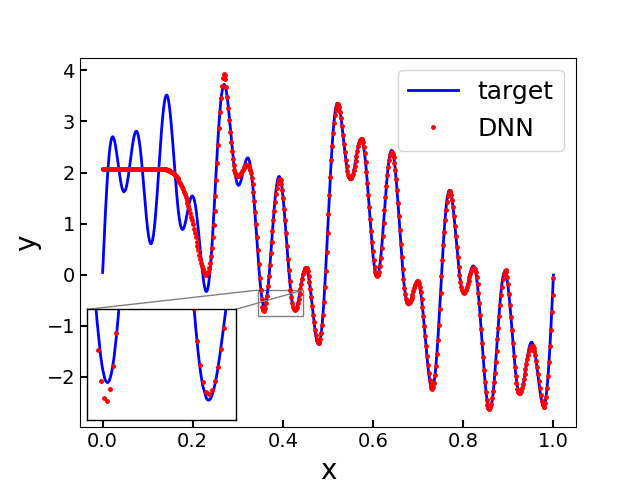} &
			\includegraphics[width=0.3\textwidth]{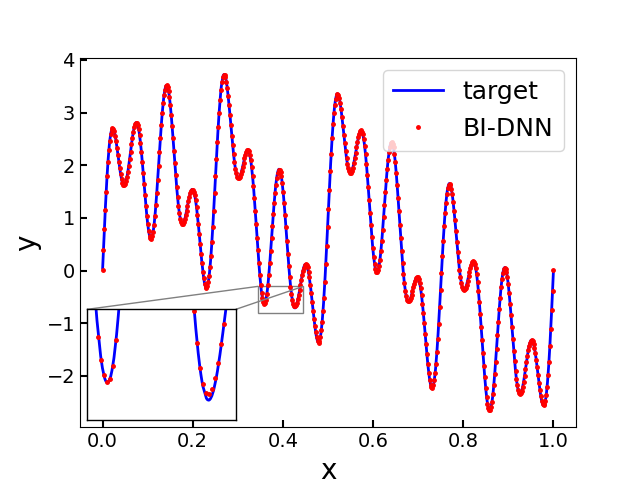} &
			\includegraphics[width=0.3\textwidth]{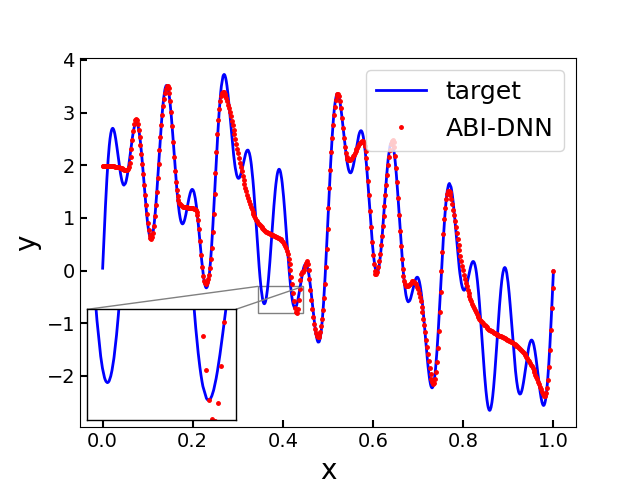} \\
			(a) DNN & (b)  BI-DNN & (c) ABI-DNN, 1st iteration \\
			\includegraphics[width=0.3\textwidth]{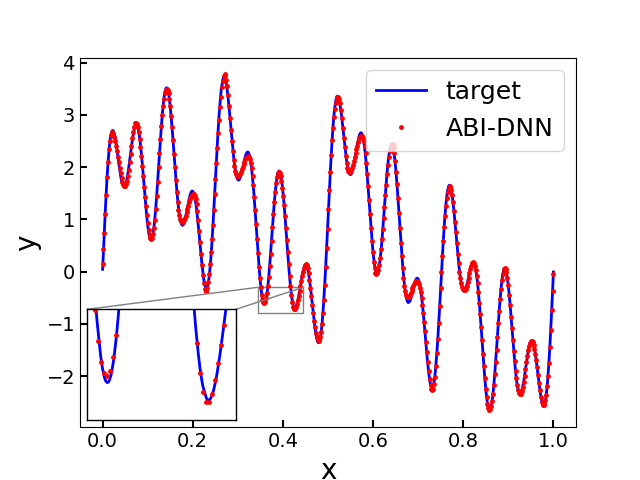} &
			\includegraphics[width=0.3\textwidth]{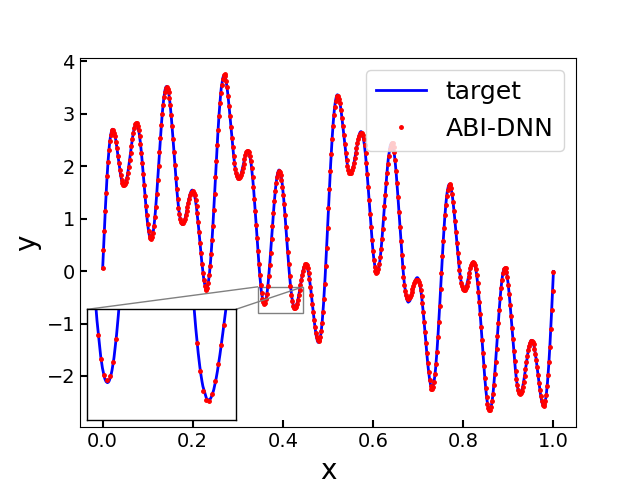} &
			\includegraphics[width=0.3\textwidth]{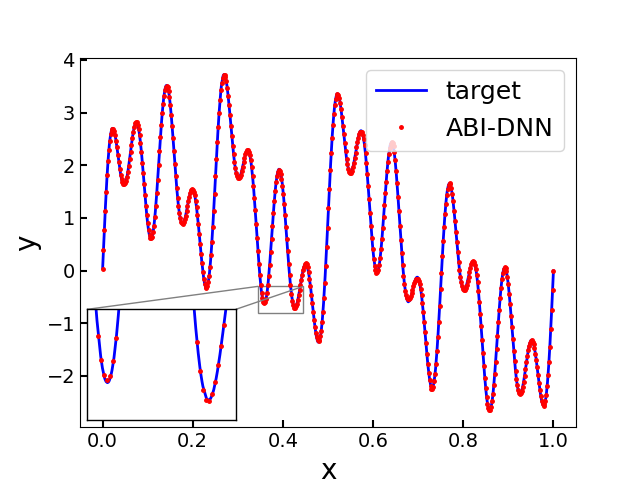} \\
			(d) ABI-DNN, 3rd iteration & (e) ABI-DNN, 5th iteration & (f) ABI-DNN, 7th iteration\\
		\end{tabular}
	\end{center}
	\caption{Illustration of the target function and approximating functions on problem \eqref{pb:1d:smooth}. Blue curves in (a-f) represent the target function, and the red dots represent approximating functions given by: (a) DNN(w=12), (b) BI-DNN(b=26), (c-f) ABI-DNN after 1st, 3rd, 5th and 7th adaptive iteration, respectively. The insets of each figure are zoomed views of approximating solutions within the region  $[0.345,0.445]\times [-0.8,-0.3]$.}
	\label{fig:fitting1d:exact:smooth}
\end{figure}
\begin{figure}[htbp]
	\begin{center}
		\begin{tabular}{ccc}
			\centering
			\includegraphics[width=0.32\textwidth]{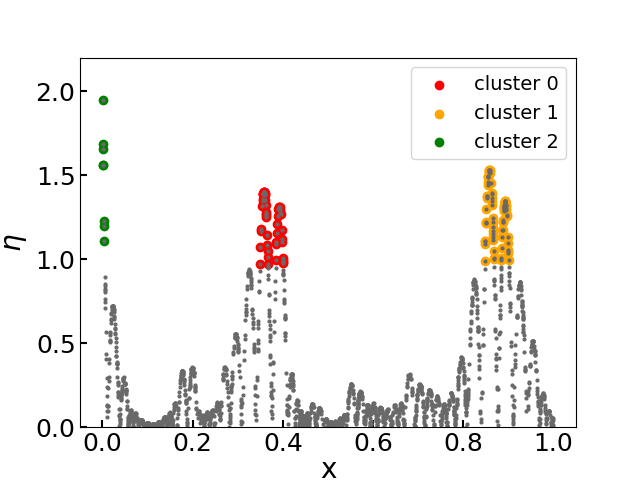}&
			\includegraphics[width=0.32\textwidth]{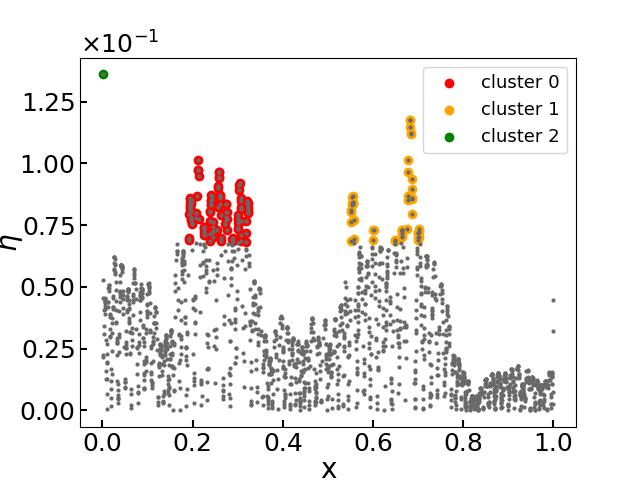}&
			\includegraphics[width=0.32\textwidth]{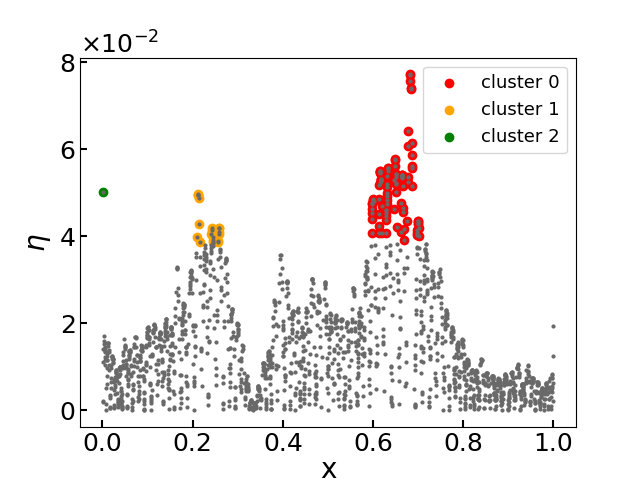}\\
			(a) after 1st iteration  & (b) after 3rd iteration   & (c) after 5th iteration
		\end{tabular}
	\end{center}
	\caption{Illustration of clustering results of ABI-DNN after 1st, 3nd, and 5th adaptive iterations on Problem \eqref{pb:1d:smooth}.}
	\label{fig:fitting1d:cluster:smooth}
\end{figure}

For the smooth target function \eqref{pb:1d:smooth}, we also start with a BI-DNN containing $10$ BI-blocks. Considering the inherent difficulty introduced by high-frequency components, we set the tolerance as $\eta_{tol}=0.015$. The adaptive process automatically terminates after $7$ iterations, yielding a final ABI-DNN model with $26$ BI-blocks. Approximation results of intermediate models are recorded in Table \ref{tab:fitting1d:ABI-DNN:smooth}.
Clearly, the ABI-DNN exhibits a consistent reduction in relative errors with increasing parameters, unlike the BI-DNN, which experiences fluctuations as shown in Fig. \ref{fig:fitting1d:params}(b).
For the purpose of comparison, the fixed BI-DNN($b=26$) and  DNN($w=12$) are trained for $70000$ epochs with their results also recorded in Table \ref{tab:fitting1d:ABI-DNN:smooth}. It can be seen that the final ABI-DNN achieves similar accuracy as its fixed counterpart with the same architecture, which is about two orders of magnitude better than the DNN with a comparable number of parameters. This demonstrates the potential of BI-DNN and ABI-DNN in dealing with complex target functions.
Fig. \ref{fig:fitting1d:exact:smooth} depicts the approximating functions provided by these three models. Specifically, from Fig. \ref{fig:fitting1d:exact:smooth}(a-b), we can see that the DNN(w=12) fails to capture the local oscillations present in $[0,0.2]$, while BI-DNN(b=26) effectively characterizes all local oscillations across the entire interval. Furthermore, the ABI-DNN exhibits an incremental learning process, progressively refining its understanding of local oscillations from a coarse scale to a more detailed fine scale, as shown in Figs.\ref{fig:fitting1d:exact:smooth}(c-f). To closely examine the adaptive enhancement processes, clustering results of ABI-DNN after the first, third, and fifth adaptive iterations are shown in Fig. \ref{fig:fitting1d:cluster:smooth}. 
Along with the marked region being enhanced, the iterative process gradually drags all pointwise error down to a smaller scale with a trend to distribute the error evenly among all training samples. 
Fig. \ref{fig:fitting1d:convergence}(b) compares the convergence processes, highlighting the rapid convergence of both BI-DNN and ABI-DNN models.

\begin{table*}[htbp]
	\centering
	\caption{Numerical results of ABI-DNN on Problem \eqref{pb:1d:smooth}. For comparison, the results of BI-DNN and DNN with a similar number of parameters are shown in the last two rows.
	}
	\begin{tabular}{c|ccccc}
		\toprule
		Model & Network structure & $\sharp$  Parameters & Testing error \\
		\midrule
		ABI-DNN(b=10) & 1-20-20-10-1 & 121 & 2.72E-01  \\
		\midrule
		ABI-DNN(b=13) & 1-26-26-13-1 & 157 & 8.59E-02  \\
		\midrule
		ABI-DNN(b=15) & 1-30-30-15-1 & 181 & 2.33E-02  \\
		\midrule
		ABI-DNN(b=18) & 1-36-36-18-1 & 217 & 1.62E-02  \\
		\midrule
		ABI-DNN(b=21) & 1-42-42-21-1 & 253 & 1.22E-02  \\
		\midrule
		ABI-DNN(b=24) & 1-48-48-24-1 & 289 & 9.47E-03  \\
		\midrule
		\textbf{ABI-DNN(b=26)} & \textbf{1-52-52-26-1} & \textbf{313} & \textbf{8.00E-03}  \\
		\midrule
		BI-DNN(b=26)  & 1-52-52-26-1 & 313 & 4.22E-03  \\
		\midrule
		DNN(w=12)          & 1-12-12-12-1 & 349 & 2.09E-01  \\
		\bottomrule
	\end{tabular}
	\label{tab:fitting1d:ABI-DNN:smooth}
\end{table*}

\subsubsection{The influence of different numbers of initial blocks in ABI-DNN}
Finally, we explore the impact of varying numbers of initial BI-blocks on the performance of ABI-DNN. The numbers of adaptive iterations executed, the numbers of parameters in the final ABI-DNN and the final relative errors achieved are summarized in Tables \ref{tab:fitting1d:blocks:singular} and \ref{tab:fitting1d:blocks:smooth} for the singular problem \eqref{pb:1d:singular} and the smooth problem \eqref{pb:1d:smooth}, respectively. To achieve the prescribed tolerance, ABI-DNNs initialized with fewer BI-blocks require more adaptive iterations than those initialized with a larger number of BI-blocks. Despite these differences, all configurations ultimately achieve similar final errors with comparable parameter counts. These findings suggest that ABI-DNN is relatively robust to the choice of the number of initial blocks. 
\begin{table}[htbp]
	\centering
	\caption{Effect of initial block counts on ABI-DNN performance for Problem \eqref{pb:1d:singular}.
	}
	\begin{tabular}{c|ccccc}
		\toprule
		\# Intial blocks& 5 & \textbf{10} & 15 & 20 \\
		\midrule
		\# Adaptive iterations & 7 & \textbf{5} & 4 & 3\\
		\midrule
		\# Parameters & 193  & \textbf{193} & 217 & 265\\ 
		\midrule
		Relative errors & 7.64E-04  & \textbf{4.51E-04} & 7.37E-04 & 6.32E-04\\  
		\bottomrule
	\end{tabular}
	\label{tab:fitting1d:blocks:singular}
\end{table}
\begin{table*}[htbp]
	\centering
	\caption{Effect of initial block counts on ABI-DNN performance for Problem \eqref{pb:1d:smooth}.
	}
	\begin{tabular}{c|cccccc}
		\toprule
		\# Intial blocks& 5 &  \textbf{10} & 15 & 20 & 25 \\
		\midrule
		\# Adaptive iterations  & 39 & \textbf{7} & 12 & 3 & 1 \\
		\midrule
		\# Parameters & 661 & \textbf{313} & 349 & 289 & 301 \\
		\midrule
		Relative errors & 9.12E-03 & \textbf{8.00E-03} & 8.87E-03 & 5.53E-03 & 8.23E-03 \\
		\bottomrule
	\end{tabular}
\label{tab:fitting1d:blocks:smooth}
\end{table*}

\subsection{Poisson Equation}
In this subsection, we apply the proposed method to problems involving singularities, e.g., point sources or geometric singularities. The following Poisson equation with a pure Dirichlet condition
\beq\label{eq:pde2d:poisson}
-\Delta u&=f,\quad \text{in}\;\;\Omega,\\
u&=g,\quad \text{on}\;\;\partial\Omega,
\eeq
is considered.
The domain $\Omega$ will be specified in each example, and the right-hand side $f$ and the boundary condition $g$ can be derived directly from the analytical solution, which will be defined later. Subsections 5.3 and 5.4 focus on employing BI-DNNs and ABI-DNNs that include two hidden layers following those containing BI-blocks. Specifically, the subnetwork $\mathbf{F}_{1:L}$, as defined in equation \eqref{eq:bi-dnn:dd}, is a three-layer fully connected neural network, where each hidden layer has the same number of neurons as the layer preceding the subnetwork.
Accordingly, the architecture of a BI-DNN is determined only by the number of BI blocks used in each dimension. With this configuration, the BI-DNN comprises a total of six layers. To ensure a fair comparison, PINNs with six layers and an equal number of neurons in each hidden layer are engaged when necessary. Within this subsection, when we compare the performance of the BI-DNN and the PINN, BI-DNNs with $5,10,15,20,25,30$ BI-blocks in each dimension and PINNs with $10,15,20,25,30,35,40,45$ neurons in each hidden layer are both trained for $50000$ epochs to observe the trend of error reduction with the increase of the number of parameters.
The training of these models is guided by the loss function outlined in equation \eqref{eq:loss}. The training data set is constructed from a random sampling of $N_{r} = 40000$ interior points and an additional $N_{b} = 100\times 4$ boundary points unless otherwise stated.
\subsubsection{Problem with one steep peak}
\begin{figure}[htbp]
	\begin{center}
		\begin{tabular}{ccc}
			\centering
			\includegraphics[width=0.32\textwidth]{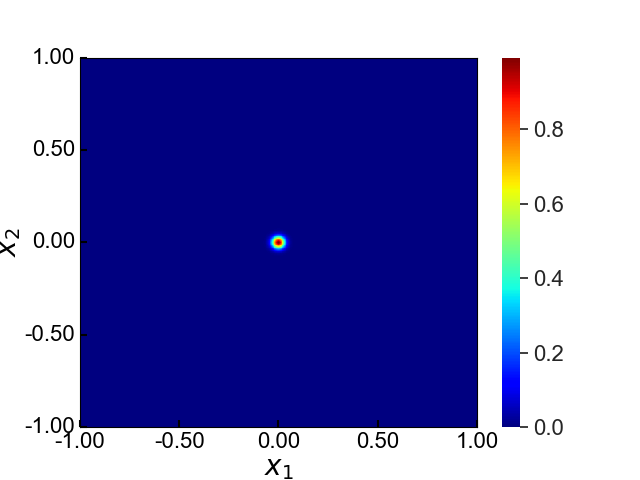} &
			\includegraphics[width=0.32\textwidth]{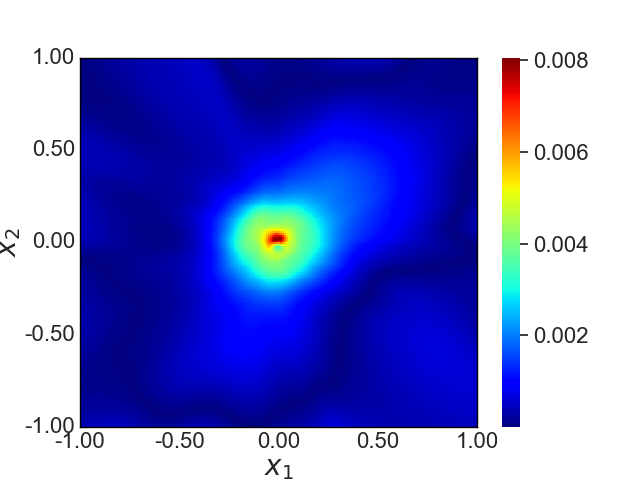} &
			\includegraphics[width=0.32\textwidth]{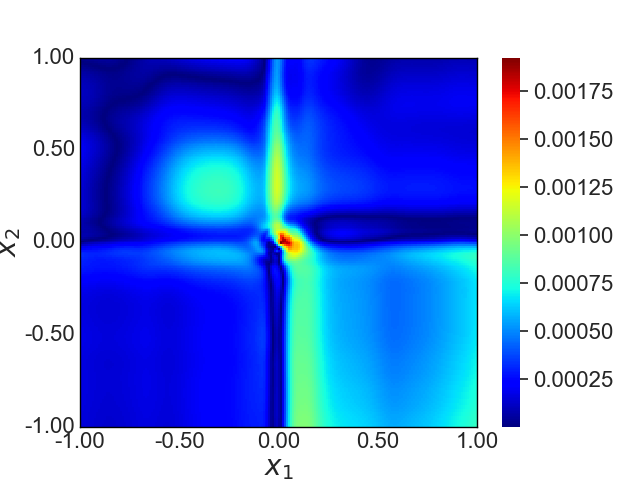} \\
			(a) Exact solution & (b) PINN & (c) BI-DNN\\
			\includegraphics[width=0.32\textwidth]{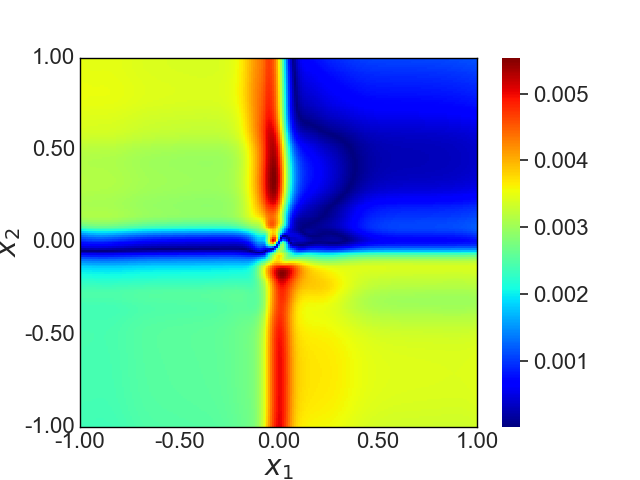} &
			\includegraphics[width=0.32\textwidth]{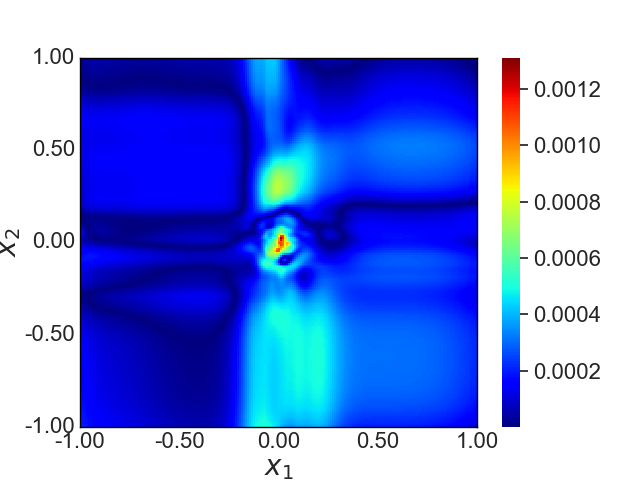} &
			\includegraphics[width=0.32\textwidth]{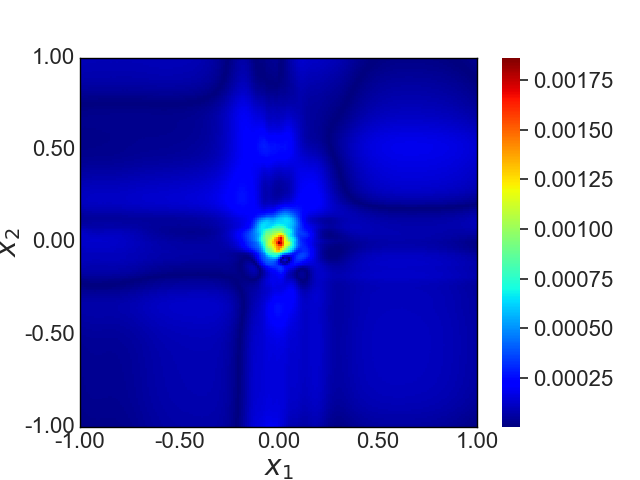} \\
			(d) ABI-DNN, 1st iteration & (e) ABI-DNN, 2nd iteration & (f) ABI-DNN, 3rd iteration
		\end{tabular}
	\end{center}
	\caption{Exact solution and pointwise absolute errors of PINN(w=19), BI-DNN(b=[12,12]) and ABI-DNN on Problem \eqref{pb:2d:onepeak}. (a) the exact solution; (b) and (c) the pointwise errors for PINN and BI-DNN, respectively; (d-f) the pointwise errors for ABI-DNN after 1st, 2nd, and 3rd adaptive iterations, respectively.
	}
	\label{fig:pde2d:pointwise:onepeak}
\end{figure}

We first consider the problem \eqref{eq:pde2d:poisson} on the quadrilateral domain $\Omega=[-1,1]^2$. The analytical solution is chosen as follows:
\be
u(x_1,x_2)= \rm e^{-1000(x_1^2+x_2^2)},
\label{pb:2d:onepeak}
\ee
which has an exponential peak at $(0,0)$ and decreases rapidly from $(0,0)$, as shown in Fig. \ref{fig:pde2d:pointwise:onepeak}(a).
\begin{figure}[htbp]
	\begin{center}
		\begin{tabular}{cc}
			\centering
			\includegraphics[width=0.45\textwidth]{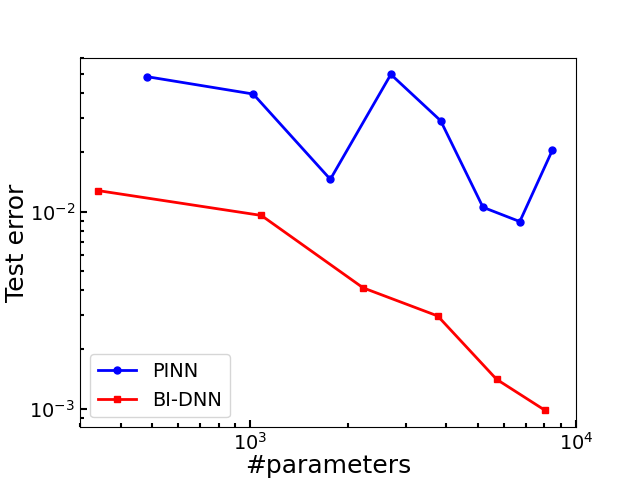}
		\end{tabular}
	\end{center}
	\caption{The change of relative $L^2$ errors of PINN and BI-DNN along with the number of trainable parameters on Problem \eqref{pb:2d:onepeak}.}
	\label{fig:pde2d:params:onepeak}
\end{figure}
\par
We first compare the performance of BI-DNN and PINN. The trend of error reduction as the number of parameters increases is depicted in Fig. \ref{fig:pde2d:params:onepeak}. With the growth in the number of parameters, BI-DNN exhibits a steady decrease in relative errors, achieving an error that is one order of magnitude lower than that of the PINN with a similar parameter count.

\begin{figure}[htbp]
	\begin{center}
		\begin{tabular}{cc}
			\centering
			\includegraphics[width=0.45\textwidth]{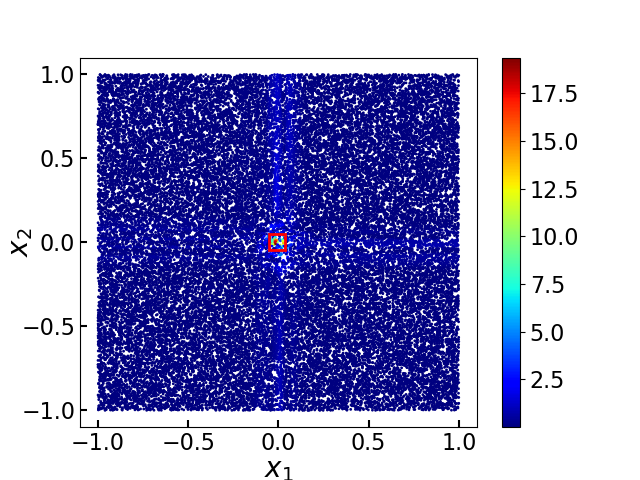}&
			\includegraphics[width=0.45\textwidth]{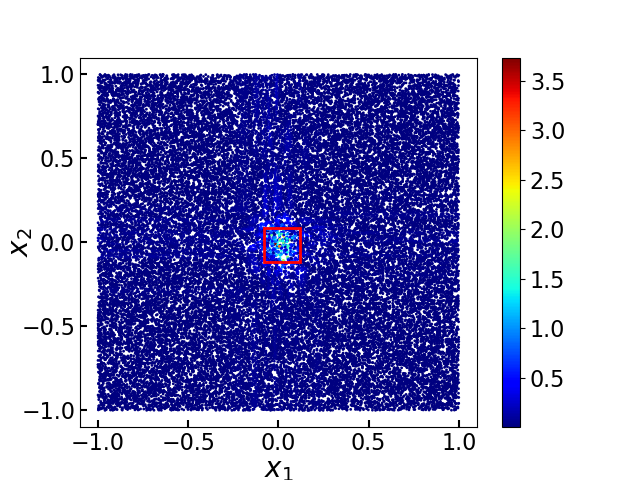}\\
			(a) after 1st iteration & (b) after 2nd iteration
		\end{tabular}
	\end{center}
	\caption{Illustration of clustering results of ABI-DNN after 1st and 2nd adaptive iterations on Problem \eqref{pb:2d:onepeak}. Each red box, centered at a cluster's centroid, has a length equal to the cluster's diameter.}
	\label{fig:pde2d:cluster:onepeak}
\end{figure}

We proceed by testing the performance of the ABI-DNN. The starting point is a BI-DNN containing $10$ BI-blocks in each dimension. For this example, we establish the tolerance at $\eta_{tol}=0.07$. The adaptive process is carried out $3$ times, with each intermediate BI-DNN being trained for $15000$ epochs, ultimately terminating at a BI-DNN with $12$ BI-blocks per dimension. The evolution of the ABI-DNN at different adaptive iterations is documented in Table \ref{tab:pde2d:ABI-DNN:onepeak}. To provide a contrast, Table \ref{tab:pde2d:ABI-DNN:onepeak} also lists the results of a fixed BI-DNN($b=[12,12]$) and a PINN($w=19$) both after $45000$ epochs of training. 
Notably, the final ABI-DNN($b=[12,12]$) achieves an error lower than the fixed BI-DNN with the same architecture and approximately an order of magnitude lower than the error yielded by the PINN with a comparable number of parameters. The results in agreement with those in Table \eqref{tab:fitting1d:ABI-DNN:singular}, further highlight the ABI-DNN's superiority in capturing sharp singularities. 
Furthermore, the pointwise absolute errors for these three models are shown in Fig. \ref{fig:pde2d:pointwise:onepeak}. As expected, testing errors are primarily concentrated around the peak, indicating the significant challenges posed by the steep gradient. The ABI-DNN demonstrates superior effectiveness in resolving sharp features compared to both the PINN and the BI-DNN. This advantage may be attributed to the ABI-DNN's accurate identification and targeted enhancement of challenging local regions, as illustrated in Fig. \ref{fig:pde2d:cluster:onepeak}.
\begin{table*}[htbp]
	\centering
	\caption{Numerical results of ABI-DNN on Problem \eqref{pb:2d:onepeak}. For comparison, the results of BI-DNN and DNN with a similar number of parameters are shown in the last two rows.
	}
	\begin{tabular}{c|ccccc}
		\toprule
		Model & Network structure & $\sharp$  Parameters & Testing error \\
		\midrule
		ABI-DNN(b=[10,10]) & 1-40-40-20-20-20-1 & 1081 & 1.88E-02  \\
		\midrule
		ABI-DNN(b=[11,11]) & 1-44-44-22-22-22-1 & 1277 & 1.02E-02  \\
		\midrule
		\textbf{ABI-DNN(b=[12,12])} & \textbf{1-48-48-24-24-24-1} & \textbf{1489} & \textbf{6.78E-03}  \\
		\midrule
		BI-DNN(b=[12,12])  & 1-48-48-24-24-24-1 & 1489 & 1.39E-02  \\
		\midrule
		PINN(w=19)          & 1-19-19-19-19-19-1 & 1597 & 5.49E-02  \\
		\bottomrule
	\end{tabular}
	\label{tab:pde2d:ABI-DNN:onepeak}
\end{table*}


\subsubsection{Problem with two steep peaks}
\begin{figure}[htbp]
	\begin{center}
		\begin{tabular}{ccc}
			\centering
			\includegraphics[width=0.32\textwidth]{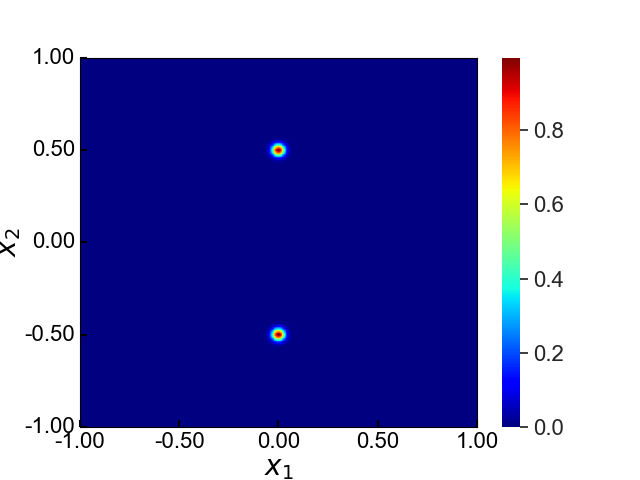} &
			\includegraphics[width=0.32\textwidth]{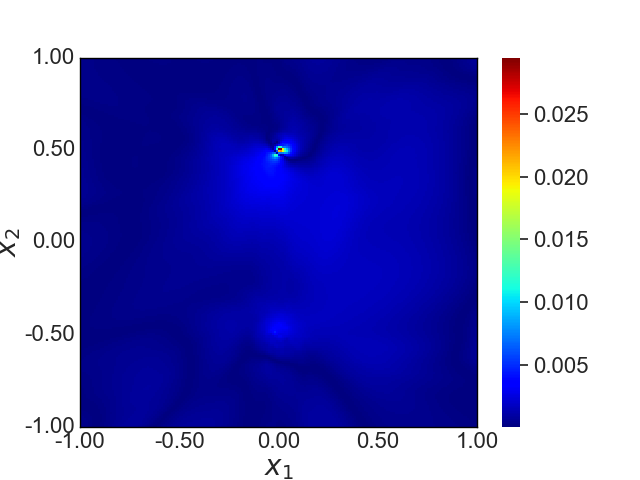} &
			\includegraphics[width=0.32\textwidth]{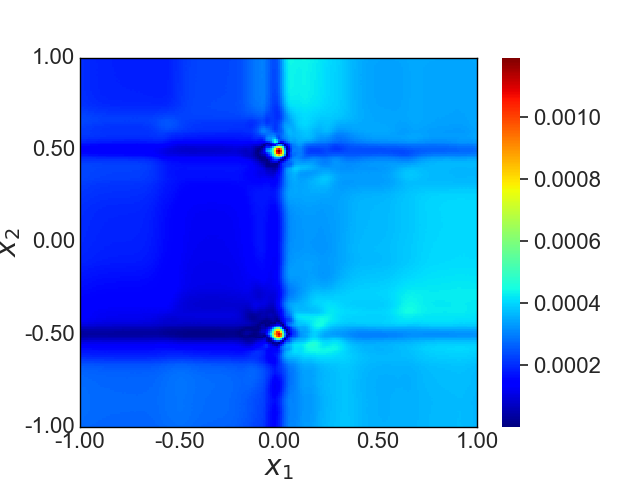} \\
			(a) Exact solution & (b) PINN & (c) BI-DNN\\
			\includegraphics[width=0.32\textwidth]{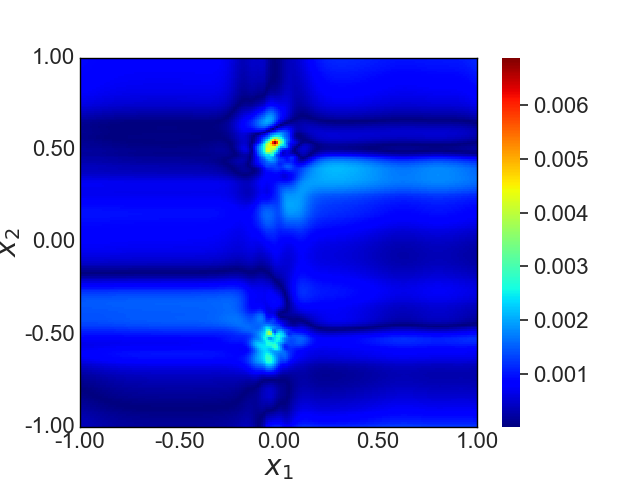} &
			\includegraphics[width=0.32\textwidth]{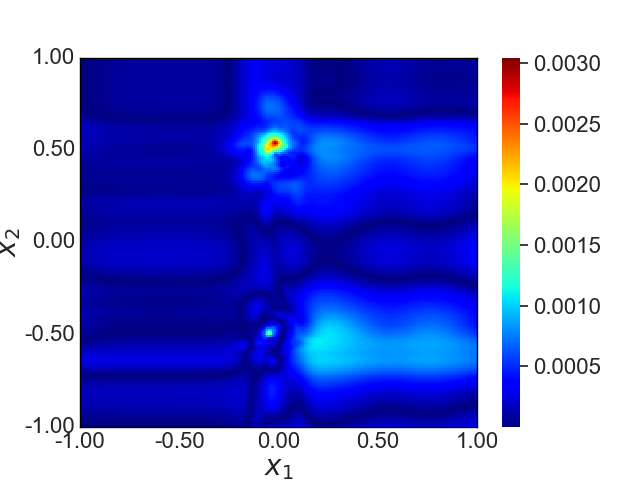} &
			\includegraphics[width=0.32\textwidth]{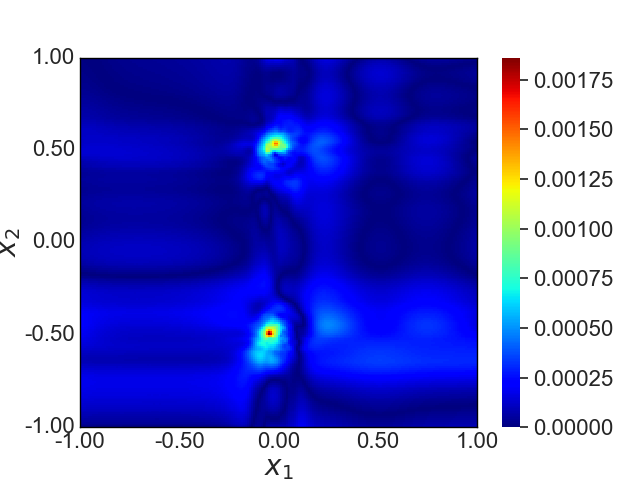}\\
			(d) ABI-DNN, 1st iteration & (e) ABI-DNN, 3rd iteration & (f) ABI-DNN, 6th iteration
		\end{tabular}
	\end{center}
	\caption{Exact solution and pointwise absolute errors of PINN(w=30), BI-DNN(b=[20,20]) and ABI-DNN on Problem \eqref{pb:2d:twopeaks}. (a) the exact solution; (b) and (c) the pointwise errors for PINN and BI-DNN, respectively; (d-f) the pointwise errors for ABI-DNN after 1st, 3rd, and 6th adaptive iterations, respectively.}
	\label{fig:pde2d:pointwise:twopeaks}
\end{figure}

To further explore the potential of ABI-DNN in dealing with multi-peak problems, we modify the exact solution to
\be
u(x_1,x_2)= \rm e^{-1000[x_1^2+(x_2-0.5)^2]}+\rm e^{-1000[x_1^2+(x_2+0.5)^2]},
\label{pb:2d:twopeaks}
\ee
which has two steep peaks at $(0,0.5)$ and $(0,-0.5)$ as shown in Fig. \ref{fig:pde2d:pointwise:twopeaks} (a).
\par
\begin{figure}[htbp]
	\begin{center}
		\begin{tabular}{c}
			\centering
			\includegraphics[width=0.5\textwidth]{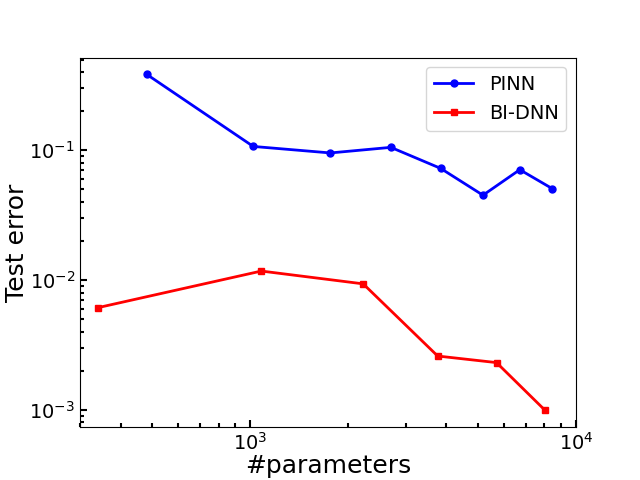}
		\end{tabular}
	\end{center}
	\caption{The change of relative $L^2$ errors of PINN and BI-DNN along with the number of trainable parameters on Problem \eqref{pb:2d:twopeaks}.}
	\label{fig:pde2d:params:twopeaks}
\end{figure}
\begin{figure}[htbp]
	\begin{center}
		\begin{tabular}{ccc}
			\centering
			\includegraphics[width=0.32\textwidth]{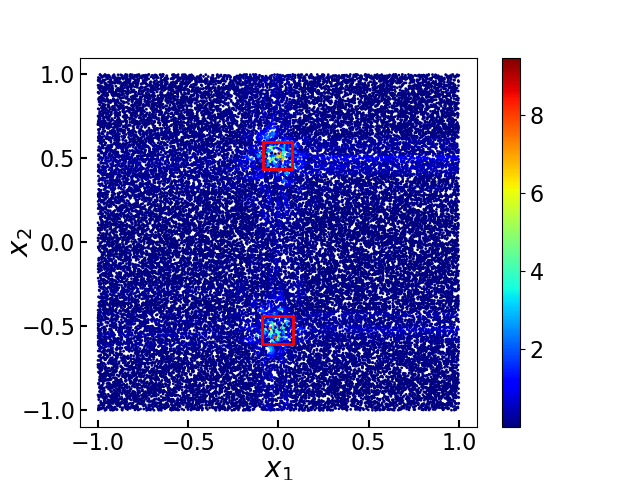}&
			\includegraphics[width=0.32\textwidth]{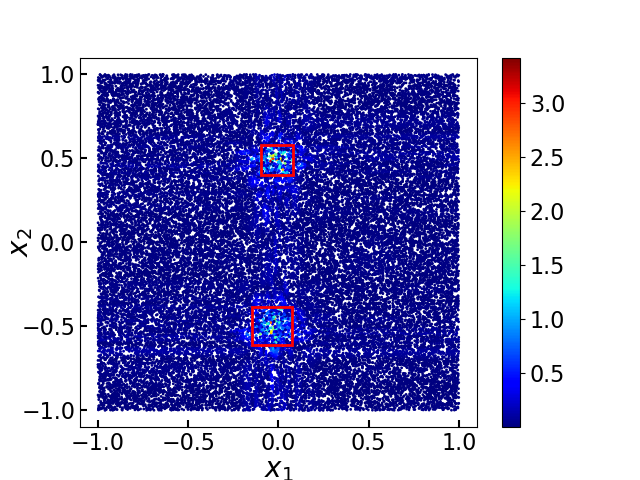}&
			\includegraphics[width=0.32\textwidth]{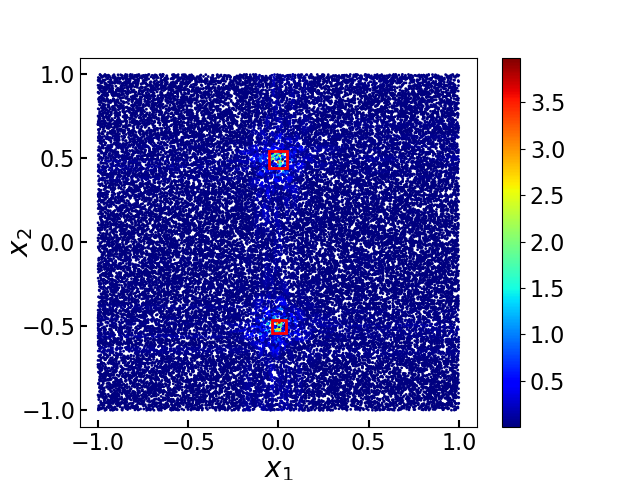}\\
			(a) after 1st iteration & (b) after 3rd iteration & (c) after 5th iteration
		\end{tabular}
	\end{center}
	\caption{Illustration of clustering results of ABI-DNN after 1st, 3rd and 5th adaptive iterations on Problem \eqref{pb:2d:twopeaks}. Each red box, centered at a cluster's centroid, has a length equal to the cluster's diameter. }
	\label{fig:pde2d:cluster:twopeaks}
\end{figure}

Fig. \ref{fig:pde2d:params:twopeaks} illustrates the error comparison between the PINN and the BI-DNN. Consistent with the trends in Fig. \ref{fig:pde2d:params:onepeak}, BI-DNNs exhibit significantly lower errors than PINNs when utilizing a similar number of parameters. 
To test the performance of ABI-DNN, we set the tolerance at $\eta_{tol}=0.1$.
The adaptive process starts with an ABI-DNN($b=[10,10]$), proceeds through $6$ iterations, and finally arrives at the ABI-DNN($b=[20,20]$), with intermediate outcomes recorded in Table \ref{tab:pde2d:ABI-DNN:twopeak}. For reference, results of the BI-DNN($b=[10,10]$) and the PINN($w=30$) after $90000$ epochs of training are also recorded. Similar to findings in the problem with one peak, both the ABI-DNN and the BI-DNN($b=[10,10]$) achieve better accuracy than the PINN($w=30$). Pointwise errors are shown in Fig. \ref{fig:pde2d:pointwise:twopeaks}, and the clustering results are presented in Fig. \ref{fig:pde2d:cluster:twopeaks}. These figures reveal that the two sharp peaks pose significant challenges for the training process, yet they are effectively spotted and enhanced during the adaptive iterations.
\begin{table*}[htbp]
	\centering
	\caption{Numerical results of ABI-DNN on Problem \eqref{pb:2d:twopeaks}. For comparison, the results of BI-DNN and DNN with a similar number of parameters are shown in the last two rows.
	}
	\begin{tabular}{c|ccccc}
		\toprule
		Model & Network structure & $\sharp$  Parameters & Testing error \\
		\midrule
		ABI-DNN(b=[10,10]) & 1-40-40-20-20-20-1 & 1081 & 2.97E-02  \\
		\midrule
		ABI-DNN(b=[12,12]) & 1-44-44-22-22-22-1 & 1489 & 1.53E-02  \\
		\midrule
		ABI-DNN(b=[14,14]) & 1-56-56-28-28-28-1 & 1961 & 1.22E-02  \\
		\midrule
		ABI-DNN(b=[16,16]) & 1-64-64-32-32-32-1 & 2497 & 8.07E-03  \\
		\midrule
		ABI-DNN(b=[18,18]) & 1-72-72-36-36-36-1 & 3097 & 6.30E-03  \\
		\midrule
		\textbf{ABI-DNN(b=[20,20])} & \textbf{1-80-80-40-40-40-1} & \textbf{3761} & \textbf{5.36E-03}  \\
		\midrule
		BI-DNN(b=[20,20])  & 1-80-80-40-40-40-1 & 3761 & 1.87E-03  \\
		\midrule
		PINN(w=30)          & 1-30-30-30-30-30-1 & 3841 & 3.92E-02  \\
		\bottomrule
	\end{tabular}
	\label{tab:pde2d:ABI-DNN:twopeak}
\end{table*}

\subsubsection{Problem with Re-entrant corner}
\begin{figure}[htbp]
	\begin{center}
		\begin{tabular}{ccc}
			\centering
			\includegraphics[width=0.32\textwidth]{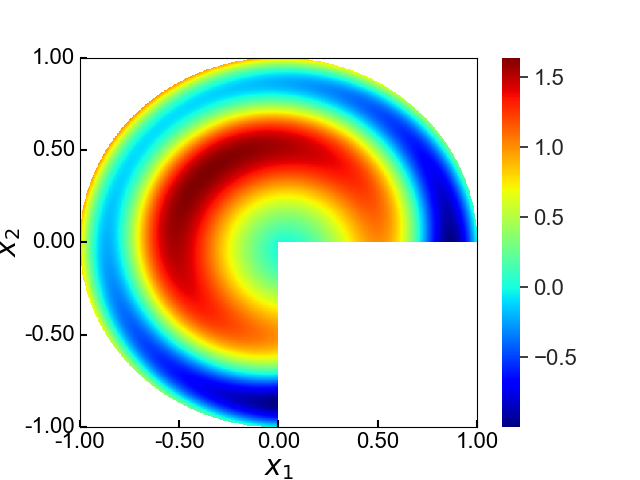}&
			\includegraphics[width=0.32\textwidth]{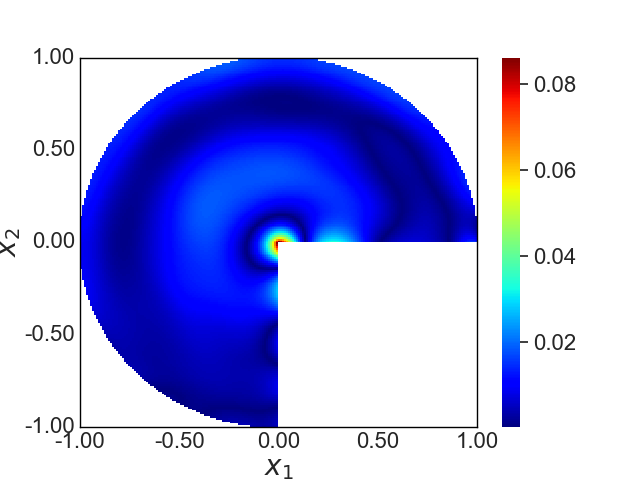} &
			\includegraphics[width=0.32\textwidth]{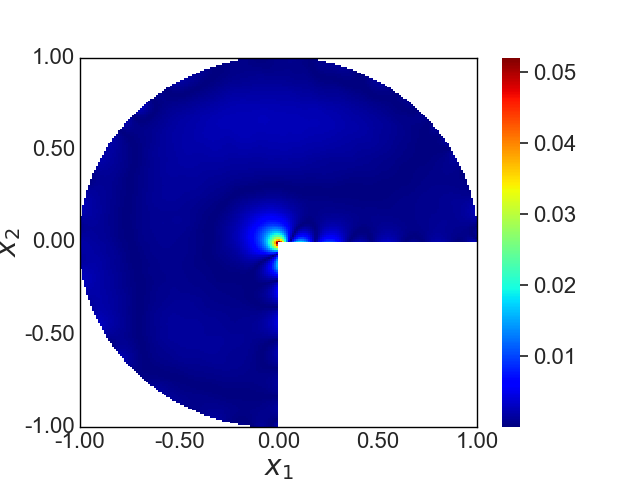}\\
			(a) Exact solution & (b) PINN & (c) BI-DNN\\
			\includegraphics[width=0.32\textwidth]{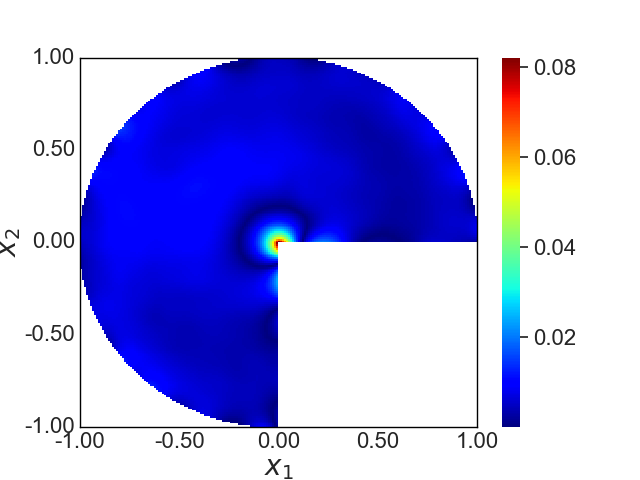} &
			\includegraphics[width=0.32\textwidth]{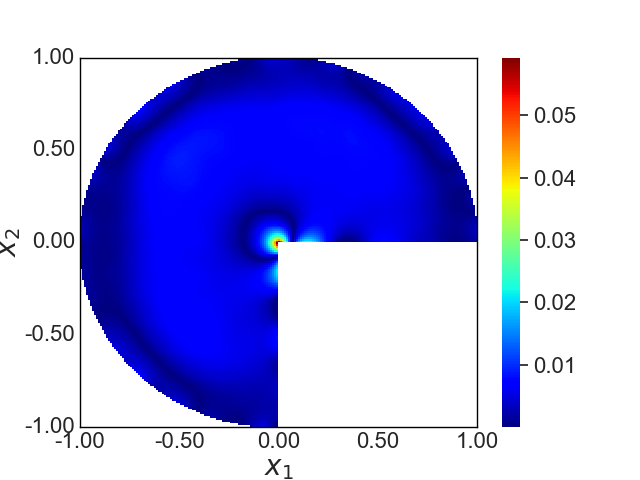} &
			\includegraphics[width=0.32\textwidth]{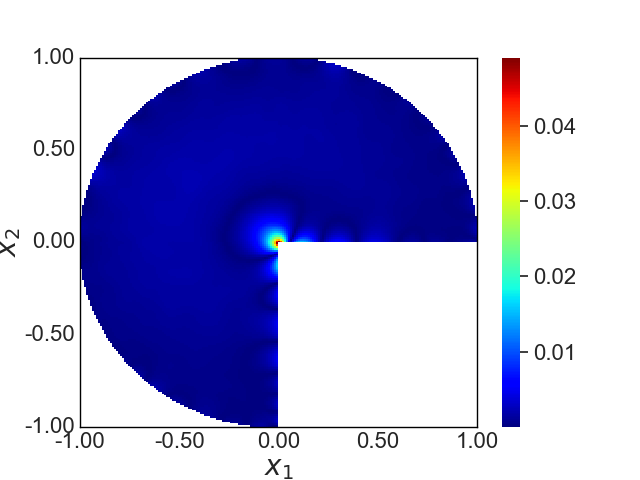} \\
			(d) ABI-DNN, 1st iteration & (e) ABI-DNN, 2nd iteration & (f) ABI-DNN, 3rd iteration
		\end{tabular}
	\end{center}
	\caption{Exact solution and pointwise absolute errors of PINN(w=21), BI-DNN(b=[14,14]) and ABI-DNN on Problem \eqref{pb:2d:Lshape}. (a) the exact solution; (b) and (c) the pointwise errors for PINN and BI-DNN, respectively; (d-f) the pointwise errors for ABI-DNN after 1st, 2nd, and 3rd adaptive iterations, respectively.
	}
	\label{fig:pde2d:pointwise:Lshape}
\end{figure}

Next, we consider the Poisson equation \eqref{eq:pde2d:poisson} defined on the domain $\Omega=\{(r,\theta)|r\in(0,1),\theta\in(0,\frac{3}{2}\pi)\}$, characterized by the presence of a re-entrant corner. The exact solution is 
\be
u(r,\theta)=r^{\frac{2}{3}}\sin(\frac{2}{3}\theta)+\sin(2\pi r^2).
\label{pb:2d:Lshape}
\ee
Note that this example without the second term $\sin(2\pi r^2)$ is usually used to test the effectiveness of the AFEM, and it represents a typical corner singularity near the re-entrant corner at the origin:  $u\in H^{\frac{5}{3}-\epsilon}$ for $\epsilon>0$. Here, we add the term $\sin(2\pi r^2)$ to make this problem more challenging. The profile of the exact solution is displayed in Fig. \ref{fig:pde2d:pointwise:Lshape}(a).

\begin{figure}[htbp]
	\begin{center}
		\begin{tabular}{c}
			\centering
			\includegraphics[width=0.5\textwidth]{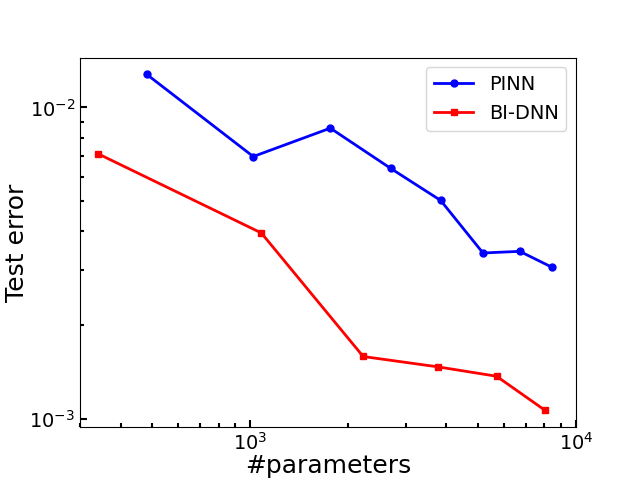}
		\end{tabular}
	\end{center}
	\caption{The change of relative $L^2$ errors of PINN and BI-DNN along with the number of trainable parameters on Problem \eqref{pb:2d:Lshape}.}
	\label{fig:pde2d:params:Lshape}
\end{figure}

\begin{figure}[htbp]
	\begin{center}
		\begin{tabular}{ccc}
			\centering
			\includegraphics[width=0.45\textwidth]{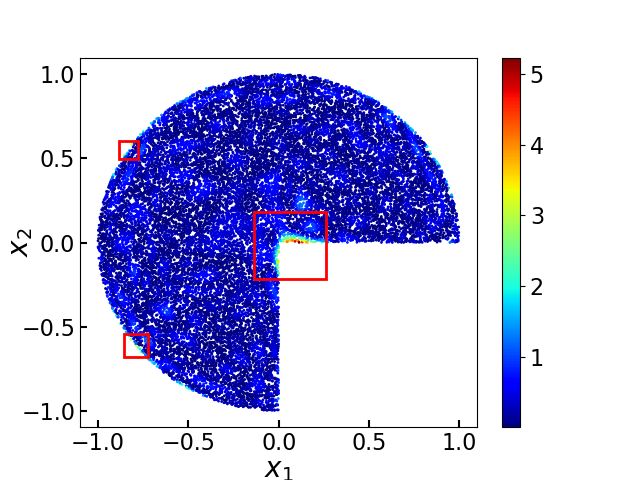}&
			\includegraphics[width=0.45\textwidth]{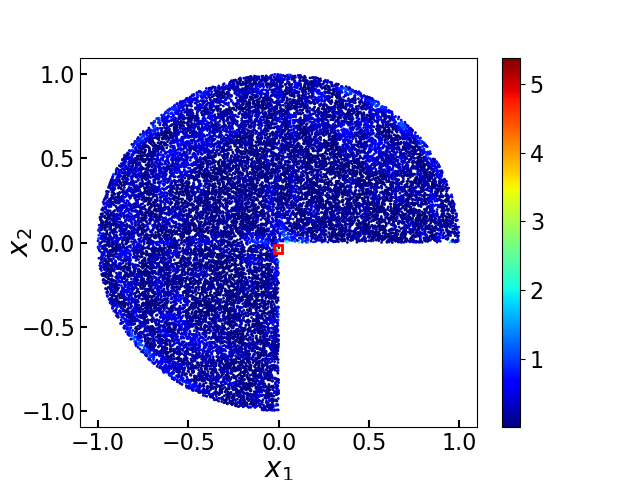}\\
			(a) after 1st iteration & (b) after 2nd iteration
		\end{tabular}
	\end{center}
	\caption{Illustration of clustering results of ABI-DNN after 1st and 2nd adaptive iterations on Problem \eqref{pb:2d:Lshape}. Each red box, centered at a cluster's centroid, has a length equal to the cluster's diameter.}
	\label{fig:pde2d:cluster:Lshape}
\end{figure}

The comparisons of BI-DNN and PINN are presented in Fig. \ref{fig:pde2d:params:Lshape}.  An obvious observation is that with the number of parameters fixed, BI-DNNs always achieve better accuracy than PINNs, which implies the architecture of the BI-DNN is more suitable for problems with corner singularities. Again, we test the effectiveness of the ABI-DNN by starting with a small model ABI-DNN($b=[10,10]$). The adaptive process experiences $4$ iterations and stops at the model ABI-DNN($b=[14,14]$) with the predefined tolerance $\eta_{tol}=0.2$ reached. The outcomes from each iterative step, as well as the results of the fixed BI-DNN($b=[14,14]$) and PINN($w=21$) through $45000$ epochs of training, are also recorded in Table \ref{tab:pde2d:ABI-DNN:Lshape}. 
As shown in Table \ref{tab:pde2d:ABI-DNN:Lshape}, the final model ABI-DNN($b=[14,14]$) achieves almost the same accuracy as the fixed model BI-DNN($b=[14,14]$), and both demonstrate approximately an order of magnitude improvement over the comparable PINN. This finding suggests that the ABI-DNN is capable of autonomously identifying a suitable network architecture by iteratively adding new BI-blocks directed by the error indicator \eqref{eq:indicator}. The pointwise absolute errors, shown in Figs. \ref{fig:pde2d:pointwise:Lshape}(b-f), highlight the singularity at the origin as the principal learning challenge. Fig. \ref{fig:pde2d:cluster:Lshape} illustrates that ABI-DNN addresses this challenge by strategically incorporating new BI-blocks around the singularity, which reveals ABI-DNN's flexibility in dynamically adjusting its architecture to adapt to the characteristics of the target function.
\begin{table*}[htbp]
	\centering
	\caption{Numerical results of ABI-DNN on Problem \eqref{pb:2d:Lshape}. For comparison, the results of BI-DNN and DNN with a similar number of parameters are shown in the last two rows.
	}
	\begin{tabular}{c|ccccc}
		\toprule
		Model & Network structure & $\sharp$  Parameters & Testing error \\
		\midrule
		ABI-DNN(b=[10,10]) & 1-40-40-20-20-20-1 & 1081 & 6.23E-03  \\
		\midrule
		ABI-DNN(b=[13,13]) & 1-52-52-26-26-26-1 & 1717 & 3.61E-03  \\
		\midrule
		\textbf{ABI-DNN(b=[14,14])} & \textbf{1-56-56-28-28-28-1} & \textbf{1961} & \textbf{2.36E-03}  \\
		\midrule
		BI-DNN(b=[14,14])  & 1-56-56-28-28-28-1 & 1961 & 2.52E-03  \\
		\midrule
		PINN(w=21)          & 1-21-21-21-21-21-1 & 1933 & 1.04E-02  \\
		\bottomrule
	\end{tabular}
	\label{tab:pde2d:ABI-DNN:Lshape}
\end{table*}

\subsection{Burgers equation}
We conclude our numerical experiments with the Burgers equation, a prominent nonlinear PDE that describes the dynamics of viscous fluids \cite{PINN, Yang_Yang_Deng_He_2024}. Owing to the presence of shockwaves, the Burgers equation presents a significant challenge to solve and remains a research focus in the field of nonlinear dynamics. 
The equation is presented below:
\begin{equation}
	\begin{aligned}
		&u_{t}+uu_{x}-(0.01/\pi)u_{xx}=0,\quad x\in[-1,1],\;t\in[0,1],\\
		&u(0,x)=-\sin(\pi x),\\
		&u(t,-1)=u(t,1)=0.
	\end{aligned}
	\label{pb:2d:burgers}
\end{equation}
It is widely recognized that despite having a smooth initial condition, the solution to the Burgers equation becomes sharp at $x = 0$ as time evolves. In this experiment, the training data set is composed of $40000$ collocation points ($N_r = 40000$) within the computational domain, $200$ boundary points ($N_b = 200$) sampled at each of the two boundaries $x = -1$ and $x = 1$, and $100$ initial points ($N_i = 100$). The testing data set comprises $256\times 100$ points as in \cite{Yang_Yang_Deng_He_2024}. We treat the time coordinate $t$ as an additional spatial dimension alongside the spatial coordinate $x$, thereby allowing the initial conditions to be considered part of the boundary conditions.

\begin{figure}[htbp]
	\begin{center}
		\begin{tabular}{c}
			\centering
			\includegraphics[width=0.45\textwidth]{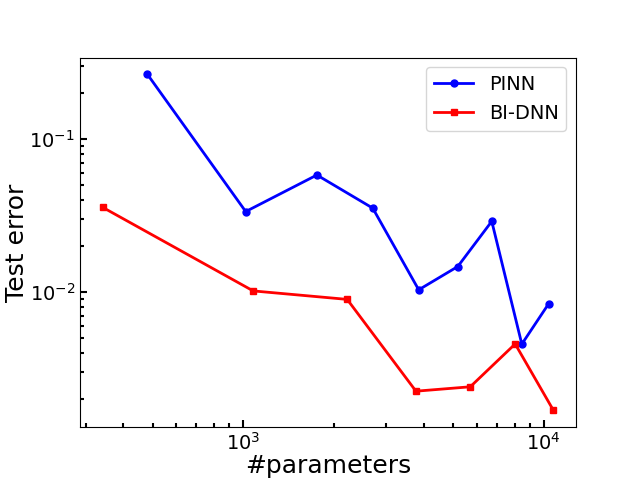}
		\end{tabular}
	\end{center}
	\caption{The change of relative $L^2$ errors of PINN and BI-DNN along with the number of trainable parameters on Problem \eqref{pb:2d:burgers}.}
	\label{fig:pde2d:params:burgers}
\end{figure}

\begin{figure}[htbp]
	\begin{center}
		\begin{tabular}{ccc}
			\centering
			\includegraphics[width=0.32\textwidth]{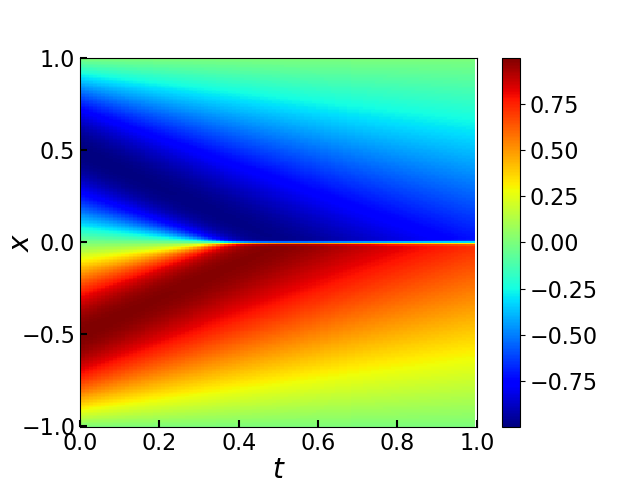}&
			\includegraphics[width=0.32\textwidth]{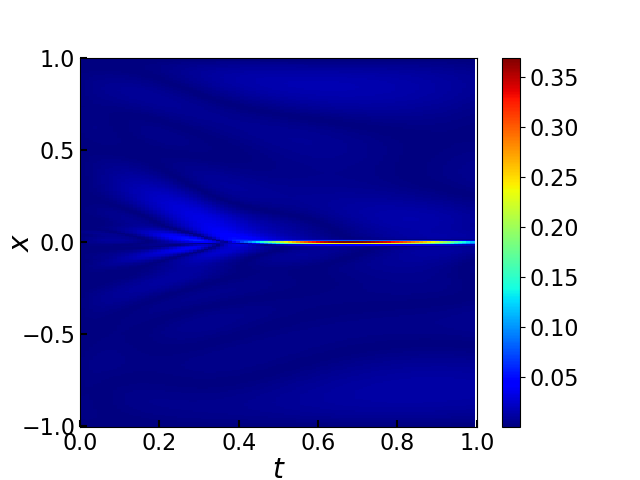} &
			\includegraphics[width=0.32\textwidth]{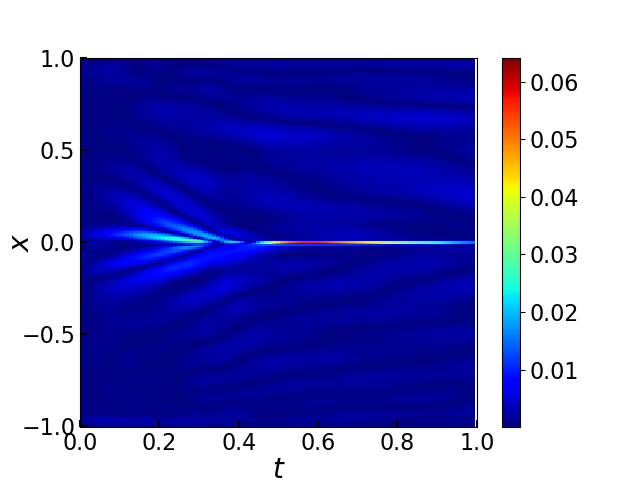}\\
			(a) Exact solution & (b) PINN & (c) BI-DNN\\
			\includegraphics[width=0.32\textwidth]{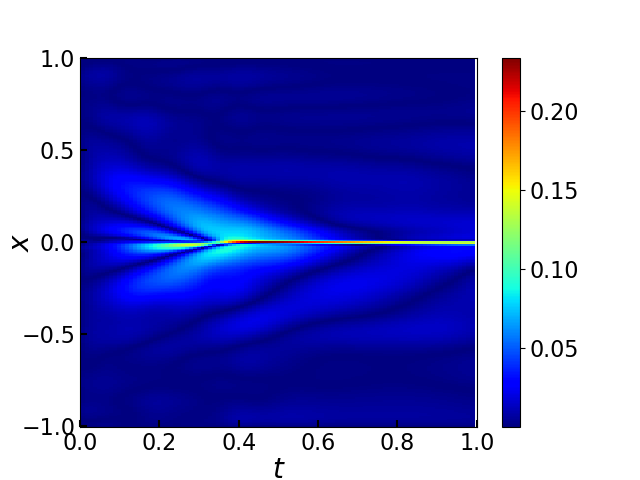} &
			\includegraphics[width=0.32\textwidth]{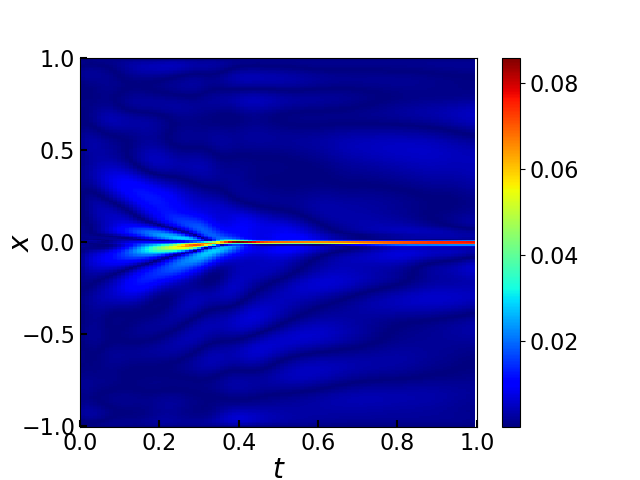} &
			\includegraphics[width=0.32\textwidth]{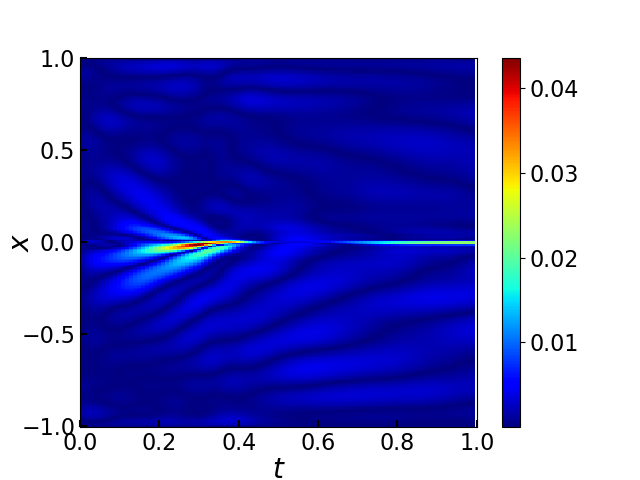} \\
			(d) ABI-DNN, 1st iteration & (e) ABI-DNN, 2nd iteration & (f) ABI-DNN, 3rd iteration
		\end{tabular}
	\end{center}
	\caption{The reference solution and pointwise absolute errors of PINN, BI-DNN, and ABI-DNN on Problem \eqref{pb:2d:burgers}. (a) the reference solution; (b) and (c) the pointwise errors for PINN(w=21) and BI-DNN(b=[14,14]), respectively; (d-f) the pointwise errors after 1st, 2nd, and 3rd adaptive iterations, respectively.
	}
	\label{fig:pde2d:pointwise:burgers}
\end{figure}

\begin{figure}[htbp]
	\begin{center}
		\begin{tabular}{ccc}
			\centering
			\includegraphics[width=0.45\textwidth]{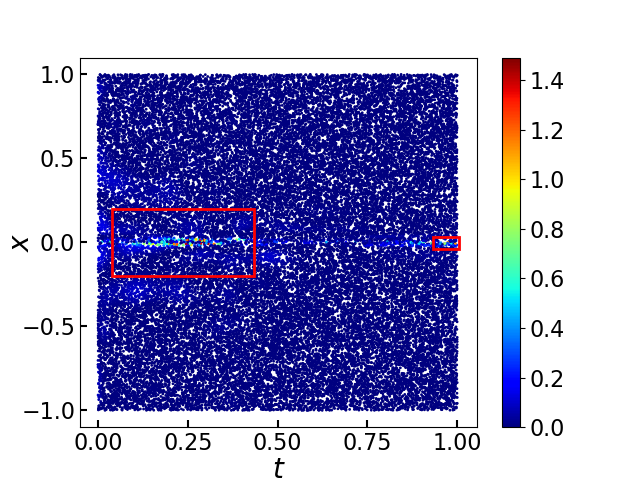}&
			\includegraphics[width=0.45\textwidth]{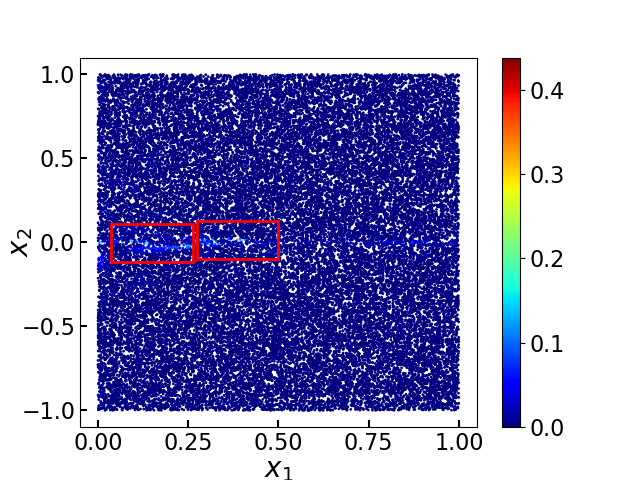}\\
			(a) after 1st iteration & (b) after 2nd iteration
		\end{tabular}
	\end{center}
	\caption{Illustration of clustering results of ABI-DNN after 1st and 2nd adaptive iterations on Problem \eqref{pb:2d:burgers}. Each red box, centered at a cluster's centroid, has a length equal to the cluster's diameter.}
	\label{fig:pde2d:cluster:burgers}
\end{figure}

We begin by comparing the numerical results of BI-DNNs with the results obtained by PINNs, all trained for $50000$ epochs. Fig. \ref{fig:pde2d:params:burgers} again shows that BI-DNNs can achieve more accurate results when a similar number of parameters is used. As in solving the Poisson equation, we also start from a small model ABI-DNN($b=[10,10$]). The adaptive enhancement process automatically stops after $3$ iterations, once the predefined tolerance $\eta_{tol}=0.005$ is achieved. Table \ref{tab:pde2d:ABI-DNN:burgers} presents a detailed comparison of the results at each iteration, alongside the outcomes for a fixed BI-DNN($b=[14,14]$) and the PINN($w=21$), each trained for $45,000$ epochs. Notably, the final ABI-DNN($b=[14,14]$) reaches the smallest error compared to the two fixed models. Fig. \ref{fig:pde2d:pointwise:burgers} shows the reference solution for the Burgers equation alongside the pointwise absolute errors of the three models from Table \ref{tab:pde2d:ABI-DNN:burgers}. These figures demonstrate that the key learning difficulty lies in the sharpness around $x=0$, and the final ABI-DNN($b=[14,14$]) captures the sharpness most precisely. 
The clustering results are shown in Fig. \ref{fig:pde2d:cluster:burgers}. Intriguingly, points highlighted by the local error indicator \eqref{eq:indicator}, although clustered around $x=0$, do not precisely overlap with the regions of high pointwise error, as depicted in Fig. \ref{fig:pde2d:pointwise:burgers}. A careful comparison between Fig. \ref{fig:pde2d:pointwise:burgers}(d) and Fig. \ref{fig:pde2d:cluster:burgers}(a), as well as Fig. \ref{fig:pde2d:pointwise:burgers}(e) and Fig. \ref{fig:pde2d:cluster:burgers}(b), reveals this discrepancy.
These observations suggest that the error indicator \eqref{eq:indicator} may not be sufficiently effective in complex problems like the Burgers equation. The development of a more sophisticated error indicator could potentially enhance the ABI-DNN's ability to produce more accurate solutions. This topic will be explored in our future studies. 

\begin{table*}[htbp]
	\centering
	\caption{Numerical results of ABI-DNN on Problem \eqref{pb:2d:burgers}.  For comparison, the results of BI-DNN and DNN with a similar number of parameters are shown in the last two rows.
	}
	\begin{tabular}{c|ccccc}
		\toprule
		Model & Network structure & $\sharp$  Parameters & Testing error \\
		\midrule
		ABI-DNN(b=[10,10]) & 1-40-40-20-20-20-1 & 1081 & 3.43E-02  \\
		\midrule
		ABI-DNN(b=[12,12]) & 1-48-48-24-24-24-1 & 1489 & 1.11E-02  \\
		\midrule
		\textbf{ABI-DNN(b=[14,14])} & \textbf{1-56-56-28-28-28-1} & \textbf{1961} & \textbf{4.89E-03}  \\
		\midrule
		BI-DNN(b=[14,14])  & 1-56-56-28-28-28-1 & 1961 & 6.18E-03  \\
		\midrule
		PINN(w=21)          & 1-21-21-21-21-21-1 & 1933 & 2.88E-02  \\
		\bottomrule
	\end{tabular}
	\label{tab:pde2d:ABI-DNN:burgers}
\end{table*}

\section{Conclusion}
In this paper, we have proposed a novel network architecture, BI-DNN, and an adaptive model, ABI-DNN, specifically designed to address problems characterized by localized features. The cornerstone of these models is the BI-block which mimics the basis function of FEM and effectively bridges the local areas that need enhancement and the neurons to be added to the network architecture. 
By leveraging BI-blocks, the ABI-DNN employs an adaptive enhancement framework that autonomously evaluates the problem's complexity and incrementally introduces new BI-blocks targeted at high-error regions, all without prior human knowledge. Consequently, the ABI-DNN can automatically generate a suitable architecture that adapts to the characteristics of the target function. A series of comprehensive numerical experiments have demonstrated the efficiency and superiority of both BI-DNN and ABI-DNN, particularly when handling problems with singularities. It is clearly observed that increased resolution is obtained in the regions of challenges. 
Future work includes developing more sophisticated error indicators, exploring the correlation between the function space determined by ABI-DNN and AFEM, and designing more efficient optimization algorithms tailored for ABI-DNN.

\section*{Acknowledgments}
This research is supported partly by National Natural Science Foundation of China with grants (12101609, 12388101), National Key R\&D Program of China with grants (2019YFA0709600, 2019YFA0709602, 2020YFA0713500), the USTC Startup Program (KY2090000141).
\bibliographystyle{unsrt}
\bibliography{ABI-DNN}
\end{document}